\documentclass[11pt]{amsart} 

\usepackage{amsmath}
\usepackage{amsfonts}
\usepackage{amssymb} 
\usepackage{amsthm} 
\usepackage{mathptmx}
\usepackage{mathrsfs}
\usepackage{latexsym}
\usepackage{times}
\usepackage[mathscr]{euscript}
\usepackage[isolatin]{inputenc}

\usepackage[pagebackref]{hyperref}
\usepackage[backrefs,alphabetic,initials]{amsrefs}

\DeclareMathAlphabet{\mathpzc}{OT1}{pzc}{m}{it}

\newtheorem{thm}{Theorem}[section]
\newtheorem{lem}[thm]{Lemma}
\newtheorem{prop}[thm]{Proposition}
\newtheorem{cor}[thm]{Corollary}

\theoremstyle{definition}
\newtheorem{defn}[thm]{Definition}
\newtheorem{ex}[thm]{Example}

\theoremstyle{remark}
\newtheorem{rem}[thm]{Remark}

\usepackage{color}

\newcommand{\slk}{\mathfrak{sl}}

\newcommand{\ok}{\mathfrak{o}}
\newcommand{\zk}{\mathfrak{z}}

\newcommand{\g}{\mathfrak{g}}

\newcommand{\hk}{\mathfrak{h}}
\newcommand{\lk}{\mathfrak{l}}

\newcommand{\sk}{\mathfrak{s}}
\newcommand{\ak}{\mathfrak{a}}

\newcommand{\qk}{\mathfrak{q}}
\newcommand{\tk}{\mathfrak{t}}
\newcommand{\bk}{\mathfrak{b}}
\newcommand{\jk}{\mathfrak{j}}
\newcommand{\uk}{\mathfrak{u}}
\newcommand{\rk}{\mathfrak{r}}

\newcommand\Bc{\mathcal{B}}
\newcommand\Lc{\mathcal L}
\newcommand\Dc{\mathcal D}

\newcommand\Ic{\mathcal I}
\newcommand\Jc{\mathcal J}
\newcommand\Pc{\mathcal P}

\newcommand\Oc{\mathcal O}

\newcommand\Nc{\mathcal N}

\newcommand{\Sb}{\mathsf{S}}

\newcommand{\Zb}{\mathsf{Z}}

\newcommand\CC{\mathbb C}
\newcommand\NN{\mathbb N}
\newcommand\ZZ{\mathbb Z}

\newcommand\PP{\mathbb P}

\newcommand{\Tr}{\operatorname{Tr}}

\newcommand{\Cs}{\mathscr C}

\newcommand{\Bs}{\mathscr B}
\newcommand{\Ds}{\mathscr D}

\newcommand{\Rs}{\mathscr R}
\newcommand{\Ss}{\mathscr S}

\newcommand{\Ns}{\mathscr N}
\newcommand{\Os}{\mathscr O}

\newcommand{\Qs}{\mathscr Q}
\newcommand{\Vs}{\mathscr V}
\newcommand{\Ws}{\mathscr W}
\newcommand{\Zs}{\mathscr Z}

\newcommand{\io}{\operatorname{\iota}}

\newcommand{\ad}{\operatorname{ad}}
\newcommand{\Ad}{\operatorname{Ad}}
\newcommand{\exx}{\operatorname{e}}

\newcommand{\im}{\operatorname{Im}} 
 
\renewcommand\dfrac{\displaystyle \frac}

\newcommand{\Id}{\operatorname{Id}} 
\newcommand{\dup}{\operatorname{dup}}
\renewcommand\hat\widehat
\renewcommand\tilde\widetilde 
\newcommand{\spa}{\operatorname{span}}
\newcommand{\GL}{\operatorname{GL}}
\newcommand{\OO}{\operatorname{O}}
\newcommand{\SO}{\operatorname{SO}}

\newcommand{\rank}{\operatorname{rank}}

\newcommand\Wedge{\bigwedge}

\newcommand\adpo{\ad_{\tt P}}

\newcommand{\perpd}{{ \ \mathop{\perp}_* \ }}
\newcommand{\oplusp}{{ \ \overset{\perp}{\mathop{\oplus}} \ }}
\newcommand{\opluspb}{{ \ \overset{\perp_B}{\mathop{\oplus}} \ }}
\newcommand{\opluspbp}{{ \ \overset{\perp_{B'}}{\mathop{\oplus}} \ }}
\newcommand{\iiso}{\overset{\mathrm{i}}{\simeq}}
\newcommand{\cb}{{\overline{C}}}
\newcommand{\cpb}{{\overline{C'}}}
\newcommand{\irm}{\mathrm{i}}

\newcommand{\diag}{\operatorname{diag}}
\newcommand{\Sso}{\Ss{}_{\mathrm{s}}}
\newcommand{\Ssoh}{\widehat{\Sso}}
\newcommand{\Ssoi}{\widehat{\Sso}^\irm}

\newcommand{\Nsi}{\widehat{\Ns}^{\mathrm{i}}}
\newcommand{\Nsh}{\widehat{{\Ns}}}
\newcommand{\Dsi}{\widehat{\Ds}^{\mathrm{i}}}
\newcommand{\Dsh}{\widehat{{\Ds}}}

\newcommand{\Si}{\Ss_{\mathrm{inv}}}
\newcommand{\Sih}{\widehat{\Si}}

\newcommand{\Ich}{\widetilde{{\Ic}}}
\newcommand{\ps}{\PP^1}
\newcommand{\dimq}{\dim_{\mathrm{q}}}
\newcommand{\ta}{{ \ \mathop{\times}\limits_{\mathrm{a}} \ }}

\begin{document}

\title[A new invariant of quadratic Lie algebras]{A new invariant of quadratic Lie algebras}

\author{Duong Minh Thanh, Georges Pinczon, Rosane Ushirobira}

\address{Institut de Mathématiques de Bourgogne, Université de
  Bourgogne, B.P. 47870, F-21078 Dijon Cedex, France}

\email{Thanh.Duong@u-bourgogne.fr}
\email{Georges.Pinczon@u-bourgogne.fr}
\email{Rosane.Ushirobira@u-bourgogne.fr}

\keywords{Quadratic Lie algebras. Invariants. Double
  extensions. Adjoint orbits. Solvable Lie algebras.}

\subjclass[2000]{17B05, 17B20, 17B30}

\date{\today}

\begin{abstract} 
  We define a new invariant of quadratic Lie algebras and give a
  complete study and classification of singular quadratic Lie
  algebras, i.e. those for which the invariant does not vanish. The
  classification is related to $\OO(n)$-adjoint orbits in $\ok(n)$.
\end{abstract}

\maketitle

\section{Introduction}

Let $\g$ be a non-Abelian quadratic Lie algebra equipped with a
bilinear form $B$. We can associate to $(\g, B)$ a canonical non-zero
3-form $I \in \Wedge^3(\g)^\g$ defined by \[ I(X,Y,Z) := B([X,Y],Z), \
\forall \ X, Y, Z \in \g.\] Let $\{ \cdot, \cdot \}$ be the
super-Poisson bracket on $\Wedge(\g)$. The 3-form $I$ satisfies (see
\cite{PU07}): \[ \{I,I\} = 0. \]

Conversely, given a quadratic vector space $(\g,B)$ and a non-zero
3-form $I \in \Wedge^3(\g)$ such that $\{I,I\} =0$, there is a
non-Abelian quadratic Lie algebra structure on $\g$ such that $I$ is
the canonical 3-form associated to $\g$ (\cite{PU07}).

Let $\Qs(n)$ be the set of non-Abelian quadratic Lie algebra
structures on the quadratic vector space $\CC^n$. We identify \[
\Qs(n) \leftrightarrow \left\{ I \in \Wedge{}^3(\CC^n) \mid \{I,I\} =
  0 \right\} \] and $\Qs(n)$ is an affine variety in $\Wedge^3(\CC^n)$
(Proposition \ref{2.6}).

The {\em $\dup$-number} of a non-Abelian quadratic Lie algebra $\g$ is
defined by \[ \dup(\g):= \dim \left( \{ \alpha \in \g^* \mid \alpha
  \wedge I = 0 \} \right), \] where $I$ is the 3-form associated to
$\g$. It measures the decomposability of the 3-form $I$ and its range
is $\{0, 1, 3\}$ (Proposition \ref{1.1}). For instance, $I$ is
decomposable if, and only if, $\dup(\g) = 3$ and the corresponding
quadratic Lie algebras are classified in \cite{PU07}, up to
i-isomorphism (i.e. isometric isomorphism). It is easy to check that
the $\dup$-number of $\g$ is invariant by i-isomorphism, that is, two
i-isomorphic quadratic Lie algebras have the same $\dup$-number (Lemma
\ref{2.1}). We shall prove in this paper, a much stronger result:

\begin{center}{\em
  the $\dup$-number of $\g$ is invariant by isomorphism.}
\end{center}

To prove this result, we need to fully understand the structure of
some particular Lie algebras. This study is interesting by itself and
we shall describe it in the sequel.

We say that a non-Abelian quadratic Lie algebra $\g$ is {\em ordinary}
if $\dup(\g) = 0$. Otherwise, $\g$ is called {\em singular}. Singular
quadratic Lie algebras are of {\em type $\Sb_1$} if their
$\dup$-number is 1 and of {\em type $\Sb_3$} if their $\dup$-number is
3.

For $n \geq 1$, let $\Os(n)$ be the set of {\em ordinary} and $\Ss(n)$
be the set of {\em singular} quadratic Lie algebra structures on
$\CC^n$. We prove the following Theorem (Propositions \ref{2.6},
\ref{2.6a} and Appendix 2):

 \bigskip

{\sf{\sc Theorem 1:}}

{\em

\begin{enumerate}

\item $\Os(n)$ is a Zariski-open subset of $\Qs(n)$.

  \smallskip

\item $\Ss(n)$ is a Zariski-closed subset of $\Qs(n)$.

  \smallskip

\item $\Qs(n) \neq \emptyset$ if, and only if, $n \geq 3$.

  \smallskip

\item $\Os(n) \neq \emptyset$ if, and only if, $n \geq 6$.

\end{enumerate}}

\bigskip

As a consequence, non-Abelian quadratic Lie algebras with dimension
higher than 6 are generically ordinary. In this work, we shall give a
complete classification of singular quadratic Lie algebras, up to
i-isomorphism and up to isomorphism.

\medskip

Let us give some details of the main results of the paper. Section
\ref{Section3} contains a preparatory study of quadratic Lie algebras
of type $\Sb_1$. It allow us to describe solvable singular Lie
algebras in terms of double extensions, a useful method introduced by
V. Kac and developed in \cite{MR} and \cite{FS}. First, we obtain
(Propositions \ref{4.4} and \ref{4.5}):

\bigskip

{\sf{\sc Theorem 2:}}

{\em

\begin{enumerate}

\item Any quadratic Lie algebra of type $\Sb_1$ is solvable and it is
  a double extension.

  \smallskip

\item A quadratic Lie algebra is singular and solvable if, and only
  if, it is a double extension.

\end{enumerate}}

\bigskip
 
What about {\bf non-solvable} singular Lie algebras? Such a Lie
algebra $\g$ can be written as
\[\g = \sk \oplusp \zk\] where $\zk$ is a central ideal of $\g$ and $\sk
\iiso \ok(3)$ equipped with a bilinear form $\lambda \kappa$ for some
non-zero $\lambda \in \CC$, where $\kappa$ is the Killing form of
$\ok(3)$ (Proposition \ref{4.5}).

In the remainder of the paper, we focus on the study of {\bf solvable
  singular Lie algebras}. We denote by $\Sso(n+2)$ the set of these
structures on $\CC^{n+2}$, by $\Ssoh(n+2)$ the set of isomorphism
classes of elements in $\Sso(n+2)$ and by $\Ssoi(n+2)$ the set of
i-isomorphism classes. Also, we denote by $\ps(\ok(n))$ the projective
space of $\ok(n)$ and by $\widetilde{\ps( \ok(n))}$ the set of orbits
of elements in $\ps(\ok(n))$ under the action induced by the
$\OO(n)$-adjoint action on $\ok(n)$. Given $\cb \in \ok(n)$, there is
an associated double extension $\g_\cb \in \Sso (n+2)$.

In Proposition \ref{4.6} and Corollary \ref{4.7}, we characterize
i-isomorphisms and isomorphisms. As a consequence, we prove the
following result, conjectured and partially proved in \cite{FS}
(Proposition 4.10):

\bigskip

{\sf{\sc Theorem 3:}}

{\em The map $\cb \to \g_\cb$ induces a bijection from
  $\widetilde{\ps(\ok(n))}$ onto $\Ssoi(n+2)$.}

\bigskip

Theorem 3 gives a remarkable relation between solvable singular
quadratic Lie algebra structures on $\CC^{n+2}$ and $\OO(n)$-adjoint
orbits in $\ok(n)$. A strong improvement to Theorem 3 will be given in
Theorem 6.

Next, we detail some particular cases. Let $\Ds(n+2)$ be the set of
{\em diagonalizable singular structures on} $\CC^{n+2}$ (i.e.  $\cb$
is a semi-simple element of $\ok(n)$) and $\Dsi(n+2)$ be the set of
i-isomorphism classes in $\Ds(n+2)$. It is clear by Theorem 3 that
$\Dsi (n+2)$ is in bijection with the well-known set of semi-simple
$\OO(n)$-orbits in $\ps(\ok(n))$ (see \cite{CM} for more details on
this set). A description of the corresponding Lie algebra structures
is given in Proposition \ref{5.5}, Corollary \ref{5.8}, Lemma
\ref{5.7} and Proposition \ref{5.10}.

Let $\Ns(n+2)$ be the set of {\em nilpotent singular structures on}
$\CC^{n+2}$, $\Nsi(n+2)$ be the set of i-isomorphism classes and
$\Nsh(n+2)$ be the set of isomorphism classes of elements in
$\Ns(n+2)$.  

In the nilpotent case, we prove (Proposition \ref{5.2}):

\bigskip

{\sf{\sc Theorem 4:}}

{\em

\begin{enumerate}

\item Let $\g$ and $\g' \in \Ns(n+2)$. Then \[\g \iiso \g' \ \text{
    if, and only if, } \ \g \simeq \g'.\] Thus $\Nsi (n+2) = \Nsh
  (n+2)$.

  \smallskip

\item Let $\widetilde{\Nc}(n)$ be the set of nilpotent $\OO(n)$-orbits
  in $\ok(n)$. Then the map $\cb \mapsto \g_\cb$ induces a bijection
  from $\widetilde{\Nc}(n)$ onto $\Nsi(n+2) = \Nsh(n+2)$.

  \smallskip

\item The set $\Nsh (n+2)$ is finite.

\end{enumerate}}

\bigskip

The classification of nilpotent $\OO(n)$-orbits in $\ok(n)$ is known
\cite{CM}. It uses deep results by Jacobson-Morosov and Kostant on
$\slk(2)$-triples in semi-simple Lie algebras. Using this
classification, we obtain a classification of $\Nsi(n+2) = \Nsh(n+2)$
in terms of {\em special} partitions of $n$ and a characterization of
the corresponding Lie algebras by means of amalgamated products of
nilpotent Jordan-type Lie algebras (Proposition \ref{5.4}).

Before working on the general case, we define the notion of an {\em
  invertible singular Lie algebra} (i.e. $\cb$ is invertible). Let
$\Si(2p+2)$ be the set of such structures on $\CC^{2p+2}$ and
$\Sih(2p+2)$ be the set of isomorphism classes of elements in
$\Si(2p+2)$. The notions of i-isomorphism and isomorphism coincide in
the invertible case as we show in Lemma \ref{5.7}.

\bigskip

Given a solvable singular Lie algebra $\g$, realized as a double
extension of $\CC^n$ by $\cb \in \ok(n)$, we consider the Fitting
components $\cb_I$ and $\cb_N$ of $\cb$ and the corresponding double
extensions $\g_I = \g_{\cb_I}$ and $\g_N = \g_{\cb_N}$ that we call
the {\em Fitting components} of $\g$. We have $\g_I$ invertible,
$\g_N$ nilpotent and we prove (Proposition \ref{6.4}):

{\sf{\sc Theorem 5:}}

{\em Let $\g$ and $\g'$ be solvable singular Lie algebras and let
  $\g_N$, $\g_I$, $\g'_N$, $\g'_I$ be their Fitting
  components. Then \[ \g \iiso \g' \ \text{ if, and only if }
  \ \begin{cases} \g_N \iiso \g'_N \\ \g_I \iiso \g'_I \end{cases}\]
  The result remains valid if we replace $\iiso$ by $\simeq \ $.

}

\bigskip

Since i-isomorphism and isomorphism are equivalent notions in the case
of nilpotent or invertible singular Lie algebras, we deduce as an
immediate Corollary:

\bigskip

{\sf{\sc Theorem 6:}}

{\em Let $\g$ and $\g'$ be solvable singular Lie algebras. Then \[ \g
  \simeq \g' \ \text{ if, and only if } \ \g \iiso \g'. \] Therefore
  $\Ssoh(n+2) = \Ssoi(n+2)$.  }

\bigskip

Theorem 6 is a really interesting and unexpected property of solvable
singular quadratic Lie algebras.

Using Theorem 5, since the study of the nilpotent case is complete, we
are left with the invertible case. First, we achieve the description
of these structures in terms of amalgamated products of Jordan-type
Lie algebras in Proposition \ref{6.5}. Then, we give a classification
of invertible $\OO(n)$-orbits in $\ok(n)$ (i.e. $\OO(n)$-orbits of
invertible elements). Let $\Ic(n)$ be the set of invertible elements
in $\ok(n)$ and $\Ich(n)$ be the set of $\OO(n)$-adjoint orbits of
elements in $\Ic(n)$. Notice that $\Ic(2p+1) = \emptyset$ (Appendix
1). Next, we consider \[\Dc = \bigcup_{r \in \NN^*} \{ (d_1, \dots,
d_r) \in \NN^r \mid d_1 \geq d_2 \geq \dots \geq d_r \geq 1 \} \] and
the map $\Phi : \Dc \to \NN$ defined by $\Phi(d_1, \dots, d_r) =
\sum_{i=1}^r d_i$. We introduce the set $\Jc_p$ of all triples
$(\Lambda, m, d)$ such that:

\begin{enumerate}

\item $\Lambda$ is a subset of $\CC \setminus \{0\}$ with
  $\sharp \Lambda \leq 2p$ and $\lambda \in \Lambda$ if, and only if,
  $-\lambda \in \Lambda$.

\smallskip

\item $m : \Lambda \to \NN^*$ satisfies $m(\lambda) = m(-\lambda)$,
  for all $\lambda \in \Lambda$ and $\sum_{\lambda \in \Lambda}
  m(\lambda) = 2p$.

\smallskip

\item $d : \Lambda \to \Dc$ satisfies $d(\lambda) = d(-\lambda)$, for
  all $\lambda \in \Lambda$ and $\Phi \circ d = m$.

\end{enumerate}

To every $\cb \in \Ic(2p)$, we can associate an element $(\Lambda, m,
d)$ of $\Jc_p$ as follows: write $\cb = S + N$ as a sum of its
semi-simple and nilpotent parts. Then $\Lambda$ is the spectrum of
$S$, $m$ is the multiplicity map on $\Lambda$ and $d$ gives the size
of the Jordan blocks of $N$. Therefore, we obtain a map $i : \Ic(2p)
\to \Jc_p$ and we prove (Proposition \ref{6.8}):

\bigskip

{\sf{\sc Theorem 7:}}

{\em The map $i : \Ic(2p) \to \Jc_p$ induces a bijection from
  $\Ich(2p)$ onto $\Jc_p$.}

\bigskip

As a Corollary, we deduce a bijection from $\Sih (2p+2)$ onto $\Jc_p /
\CC^*$ (Proposition \ref{6.11}) where the action of $\mu \in \CC^* = \CC
\setminus \{0\}$ on $\Jc_p$ is defined by \[ \mu \cdot (\Lambda, m, d)
:= (\mu \Lambda, m', d'), \ \text{ with } m'(\mu \lambda) = m
(\lambda) \ \text{ and } \ d'( \mu \lambda) = d(\lambda), \ \forall \
\lambda \in \Lambda \]

Combine Theorems 5, 4 and 7 to obtain a complete classification of
$\Ssoi(n) = \Ssoh(n)$. As a by-product, we also obtain a complete
classification of $\OO(n)$-orbits in $\ok(n)$, a result which is
certainly known, but for which we have no available reference.

Finally, as a consequence of the preceding results, we prove in
Section 7 (Proposition \ref{7.3}):

\bigskip

{\sf{\sc Theorem 8:}}

{\em The $\dup$-number is invariant under isomorphism, i.e. if \[\g
  \simeq \g' \ \text { then } \ \dup(\g) = \dup(\g').\]}

\bigskip

This result is rather unexpected. It is obtained through a computation
of centromorphisms in the reduced singular case (Proposition
\ref{7.2}).

We also obtain the quadratic dimension of $\g$ \cite{BB} in this
case: \[ \dimq (\g) = 1 + \dfrac{\dim(\Zs(\g))
  (1+\dim(\Zs(\g))}{2}, \] where $\Zs(\g)$ is the center of $\g$.

There are two Appendix. In the first one, we collect some well-known
useful properties of elements of $\ok(n)$, shorts proofs are given for
the sake of completeness. In Appendix 2, we show that $\Os(5) =
\emptyset$ and describe $\Qs(5)$ up to i-isomorphism.

\section{Preliminaries} \label{prelim}

\subsection{} All vector spaces considered in the paper are
finite-dimensional complex vector spaces.

Given a vector space $V$, we denote by $V^*$ its dual space. Given a
subset $X$ of $V$, $X^\perpd$ denotes the {\em orthogonal subspace} of
$X$ in $V^*$.

We denote by $\Lc(V)$ the {\em algebra of linear operators} of $V$, by
$\GL(V)$ the {\em group of invertible operators} in $\Lc(V)$, by
${}^tA$ the {\em transpose} of an operator $A \in \Lc(V)$ and by
$\Wedge(V)$ the ($\ZZ$-graded) {\em Grassmann algebra} of
skew-symmetric multilinear forms on $V$, i.e. $\Wedge(V)$ is the
exterior algebra of $V^*$. Recall that given an isomorphism $A$
between two vector spaces $V$ and $V'$, there is an algebra
isomorphism from $\Wedge(V')$ onto $\Wedge(V)$ that extends the
transpose ${}^tA : V'^* \to V^*$ and that we also denote by ${}^tA$.

\subsection{} Let $I \in \Wedge^k(V)$, for $k \geq 1$. We introduce
two subspaces of $V^*$:
\begin{eqnarray*}
  \Vs_I &:=& \{\alpha \in V^* \mid \alpha \wedge I = 0 \}\\
  \Ws_I &:=& \{v \in V \mid \io_v(I) = 0\}^\perpd = \{ \io_{v \wedge v'} (I) \mid v, v' \in V\}
\end{eqnarray*}
where $\io_v$ is the derivation of $\Wedge(V)$ defined by:
\[\io_v(\Omega)(v_1, \dots, v_{r-1}) = \Omega(v, v_1, \dots, v_{r-1}),
\forall \ \Omega \in \Wedge{}^r(V), v_1, \dots, v_{r-1} \in V.\]

The following result is well known, see for instance \cite{Bour}.

\begin{prop} \label{1.1} Let $I \in \Wedge^k(V)$, $I \neq 0$. Then:

\begin{enumerate}
\item $\Vs_I \subset \Ws_I$, $\dim(\Vs_I) \leq k$ and $\dim(\Ws_I) \geq k$.

\item If $\{\alpha_1, \dots, \alpha_r \}$ is a basis of $\Vs_I$, then
  $\alpha_1 \wedge \dots \wedge \alpha_r$ divides $I$. Moreover, $I$
  belongs to the $k$-th exterior power of $\Ws_I$, also denoted by
  $\Wedge^k(\Ws_I)$.

\item $I$ is decomposable if, and only if, $\dim(\Vs_I) = k$ or
  $\dim(\Ws_I) = k$. In this case, $\Vs_I = \Ws_I$ and if $\{\alpha_1,
  \dots, \alpha_k \}$ is a basis of $\Vs_I$, one has for some non-zero
  $\lambda \in \CC$, \[ I = \lambda \alpha_1 \wedge \dots \wedge
  \alpha_k.\]

\end{enumerate}

\end{prop}

\subsection{} A vector space $V$ equipped with a non-degenerate
symmetric bilinear form $B$ is called a {\em quadratic vector
  space}. In this case, there is an isomorphism $\phi$ from $V$ onto
$V^*$ defined by \[\phi(v)(v') := B(v,v'), \ \forall \ v, v' \in V.\]
Given a subspace $W$ of $V$, we denote by $W^\perp$ the {\em
  orthogonal subspace} of $W$ in $V$ with respect to the bilinear form
$B$. One has $V = W \oplus W^\perp$ if, and only if, the restriction
$B|_{W \times W}$ is non degenerate and in this case, we use the
notation \[ V = W \oplusp W^\perp.\]

\subsection{} \label{1.3} Let $(V, B)$ and $(V', B')$ be two quadratic
vector spaces. An {\em isometry} is a bijective map $A : V \to V'$
that satisfies \[B'(A(v), A(w)) = B(v,w), \ \forall \ v, w \in V.\] We
denote by $A^* \in \Lc(V)$ the {\em adjoint map} of an element $A \in
\Lc(V)$. Remark that $A$ is an isometry of $V$ if, and only if,
$A^{-1} = A^*$.

The {\em group of isometries} of $V$ is denoted by $\OO(V, B)$ (or
simply $\OO(V)$) and its Lie algebra is denoted by $\ok(V,B)$ (or
simply $\ok(V)$). An element $A$ of $\ok(V) \subset \Lc(V)$ satisfies
$A^* = -A$ (that means $A$ is skew-symmetric with respect to
$B$). Notice that $\Tr(A) = 0$, for all $A \in \ok(V)$. The {\em
  adjoint action} $\Ad$ of $\OO(V)$ on $\ok(V)$ is given by \[\Ad_U(C)
:= U C U^{-1}, \ \forall \ U \in \OO(V), C \in \ok(V).\] We denote by
$\Oc_C$, the {\em orbit} of an element $C \in \ok(V)$.

Let $V = \CC^n$. Consider the canonical basis $\Bc = \{ E_1, \dots,
E_n \}$ of $V$. If $n$ even, $n = 2p$, write $\Bc = \{E_1, \dots, E_p,
F_1, \dots, F_p\}$ and if $n$ is odd, $n = 2p+1$, write $\Bc = \{E_1,
\dots, E_p, G, F_1, \dots, F_p \}$. The {\em canonical bilinear form}
$B$ on $V$ is defined by:
\begin{itemize}

\item if $n = 2p$:
  \[ B(E_i, F_j) = \delta_{ij}, B(E_i, E_j) = B(F_i, F_j) = 0, \
  \forall \ 1 \leq i,j \leq p \]

\item if $n = 2p+1$:
  \[ \begin{cases} B(E_i, F_j) = \delta_{ij}, B(E_i, E_j) = B(F_i,
    F_j) = 0, \ \forall \ 1 \leq i,j \leq p \\ B(E_i, G) = B(F_j, G) =
    0, \\ B(G,G) = 1\end{cases} \]

\end{itemize}

In that case, $\OO(n)$ stands for $\OO(\CC^n, B)$ and $\ok(n)$ stands
for $\ok(\CC^n, B)$.

Finally, if $V$ is an $n$-dimensional quadratic vector space, then $V$
is isometrically isomorphic (i-isomorphic) to the quadratic space
$\CC^n$ \cite{Bour59}.

\subsection{} \label{1.4} Let $(V,B)$ be a quadratic vector space. We
define the super-Poisson bracket on $\Wedge(V)$ as follows (see
\cite{PU07} for details): fix an orthonormal basis $\{v_1, \dots,
v_n\}$ of $V$. Then
\[ \{\Omega, \Omega '\} := (-1)^{k+1} \sum_{j=1}^n \io_{v_j}(\Omega)
\wedge \io_{v_j}(\Omega'), \ \forall \ \Omega \in \Wedge{}^k(V), \Omega'
\in \Wedge(V).\]

For instance, if $\alpha \in V^*$, one has \[\{\alpha, \Omega\} =
\io_{\phi^{-1}(\alpha)}(\Omega), \ \forall \ \Omega \in \Wedge(V),\]
and if $\alpha' \in V^*$, $\{\alpha, \alpha' \} = B(\phi^{-1}(\alpha),
\phi^{-1}(\alpha'))$.  This definition does not depend on the choice
of the basis.

For any $\Omega \in \Wedge^k(V)$, define $\adpo(\Omega)$ by \[
\adpo(\Omega)\left(\Omega'\right) : = \{ \Omega, \Omega'\}, \ \forall
\ \Omega' \in \Wedge(V).\] Then $\adpo(\Omega)$ is a super-derivation
of degree $k-2$ of the exterior algebra \linebreak $\Wedge(V)$. One
has: \[\adpo(\Omega)\left(\{\Omega', \Omega''\}\right) = \{
\adpo(\Omega)(\Omega'), \Omega''\} + (-1)^{kk'} \{ \Omega',
\adpo(\Omega)(\Omega'')\}, \] for all $\Omega' \in \Wedge^{k'}(V)$,
$\Omega '' \in \Wedge(V)$. That implies that $\Wedge(V)$ is a graded
Lie algebra for the super-Poisson bracket.

\subsection{} \label{1.5} A {\em quadratic Lie algebra} $(\g,B)$ is a
quadratic vector space $\g$ equipped with a bilinear form $B$ and a
Lie algebra structure on $\g$ such that $B$ is invariant (that means,
$B([X,Y], Z) = B(X, [Y,Z])$, for all $X$, $Y$, $Z \in \g$).

If $(\g,B)$ is a quadratic Lie algebra, recall that \[[\g,\g] =
\Zs(\g)^\perp\] where $\Zs(\g)$ is the center of $\g$. There is a
canonical invariant $I \in \Wedge^3(\g)$ defined by \[I(X,
Y,Z):=B([X,Y],Z), \ \forall \ X, Y, Z \in \g.\] This invariant
satisfies $\{I, I \} =0$ (see \cite{PU07}) and it is easy to check
that \[ \Ws_I = \phi\left( [\g, \g ] \right).\] We say that $I$ is {\em
  the 3-form associated to } $\g$.

On the other hand, given a quadratic vector space $(\g, B )$ and $I
\in \Wedge^3(\g)$, define \[ [X, Y] := \phi^{-1} \left( \io_{X \wedge
    Y} (I) \right), \ \forall \ X, Y \in \g.\] This bracket satisfies
the Jacobi identity if, and only if, $\{I, I \} = 0$ \cite{PU07}. In
this case, $\g$ becomes a quadratic Lie algebra with invariant
bilinear form $B$.

\begin{defn}
  Let $(\g,B)$ and $(\g',B')$ be two quadratic Lie algebras.  We say
  that $(\g,B)$ and $(\g',B')$ are {\em isometrically isomorphic} (or
  {\em i-isomorphic}) if there exists a Lie algebra isomorphism $A$
  from $\g$ onto $\g'$ satisfying \[ B'(A(X), A(Y)) = B(X,Y), \
  \forall \ X, Y \in \g.\] In other words, $A$ is an i-isomorphism if
  it is a Lie algebra isomorphism and an isometry. We write $\g \iiso
  \g'$.
\end{defn}

Consider two quadratic Lie algebras $(\g,B)$ and $(\g,B')$ (same Lie
algebra) with $B' = \lambda B$, $\lambda \in \CC$, $\lambda \neq
0$. They are not necessarily i-isomorphic, as shown by the example
below:

\begin{ex}
  Let $\g = \ok(3)$ and $B$ its Killing form. Then $A$ is a Lie
  algebra automorphism of $\g$ if, and only if, $A \in \OO(\g)$. So
  $(\g,B)$ and $(\g,\lambda B)$ cannot be i-isomorphic if $\lambda
  \neq 1$.
\end{ex}

\section{The dup number of a quadratic Lie algebra}

\subsection{} Let $\g$ and $\g'$ be quadratic Lie algebras with
associated invariants $I$ and $I'$ (see (\ref{1.5})). The following
Lemma is straightforward:

\begin{lem} \label{2.1} Let $A$ be an i-isomorphism from $\g$ onto
  $\g'$. Then $I = {}^tA (I')$, $\Vs_I = {}^tA(\Vs_{I'})$ and $\Ws_I
  ={}^tA(\Ws_{I'})$.
\end{lem}

It results from the previous Lemma that $\dim(\Vs_I)$ and
$\dim(\Ws_I)$ are invariant under i-isomorphisms. This is not new for
$\dim(\Ws_I)$, since $\dim(\Ws_I) = \dim \left( [\g, \g] \right)$.

For $\dim(\Vs_I)$, to our knowledge this fact was not remarked up to
now, so we introduce the following definition:

\begin{defn} Let $\g$ be a quadratic Lie algebra. The {\em $\dup$
    number} $\dup(\g)$ is defined by \[ \dup(\g) := \dim(\Vs_I).\]
\end{defn}

\begin{rem} \label{2.3} By Proposition \ref{1.1}, when $\g$ is
  non-Abelian, one has $\dup(\g) \leq 3$. Actually $\dup(\g) \in \{0,
  1, 3 \}$. Notice that $\dim(\Ws_I) \geq 3$, so $\dim\left([\g,\g]
  \right) \geq 3$ (see \cite{PU07}), a simple but rather interesting
  remark.
\end{rem}

\subsection{} \label{2.2s} We shall use the decomposition result
below:

\begin{prop} \label{2.8} \cite{PU07} \hfill

  Let $(\g,B)$ be a non-Abelian quadratic Lie algebra. Then there
  exists a central ideal $\zk$ and an ideal $\lk \neq \{0\}$ such
  that:

\begin{enumerate}

\item $\g = \zk \oplusp \lk$

\item $\left( \zk, B|_{\zk \times \zk} \right)$ and $\left(\lk,
    B|_{\lk \times \lk} \right)$ are quadratic Lie algebras. Moreover,
  $\lk$ is non-Abelian.

\item The center $\Zs(\lk)$ is totally isotropic, i.e. $\Zs(\lk)
  \subset [\lk, \lk]$.

\item Let $\g'$ be a quadratic Lie algebra and $A : \g \to \g'$ be a
  Lie algebra isomorphism. Then \[ \g' = \zk' \oplusp \lk'\] where
  $\zk' = A(\zk)$ is central, $\lk' = A(\zk)^\perp$, $\Zs(\lk')$ is
  totally isotropic and $\lk$ and $\lk'$ are isomorphic. Moreover if
  $A$ is an i-isomorphism, then $\lk$ and $\lk'$ are i-isomorphic.

\end{enumerate}

\end{prop}

\begin{proof}
  We prove (4) : recall that $\zk$ is any complementary subspace of
  $\Zs(\g) \cap [\g,\g]$ in $\Zs(\g)$ (see \cite{PU07}) and that $\lk$
  is defined as the orthogonal subspace of $\zk$, $\lk = \zk^\perp$.

  One has $A(\Zs(\g) \cap [\g, \g]) = \Zs(\g') \cap [\g', \g']$ and
  $\Zs(\g') = \zk' \oplus (\Zs(\g') \cap [\g', \g'])$. Therefore
  $\lk'$ satisfies $\g' = \zk' \oplusp \lk'$ and $\Zs(\lk')$ is
  totally isotropic. Since $A$ is an isomorphism from $\zk$ onto
  $\zk'$, $A$ induces an isomorphism from $\g / \zk$ onto $\g' /
  \zk'$, and it results that $\lk$ and $\lk'$ are isomorphic Lie
  algebras. Same reasoning works for $A$ i-isomorphism.
\end{proof}

It is clear that $\zk = \{0\}$ if, and only if, $\Zs(\g)$ is totally
isotropic and that \[\dup(\g) = \dup(\lk).\]

\begin{defn} A quadratic Lie algebra $\g$ is {\em reduced} if:

\begin{enumerate}

\item $\g \neq \{0\}$

\item $\Zs(\g)$ is totally isotropic.

\end{enumerate}

\end{defn}

Notice that a reduced quadratic Lie algebra is necessarily
non-Abelian.

\subsection{} We separate non-Abelian quadratic Lie algebras as follows:

\begin{defn} \hfill

Let $\g$ be a non-Abelian quadratic Lie algebra. 

\begin{enumerate}

\item $\g$ is an {\em ordinary} quadratic Lie algebra if $\dup(\g) =
  0$.

\item $\g$ is a {\em singular} quadratic Lie algebra if $\dup(\g) \geq
  1$. 

\begin{itemize}

\item[(i)] $\g$ is a {\em singular} quadratic Lie algebra of {\em type
    $\Sb_1$} if $\dup(\g) = 1$.

\item[(i)] $\g$ is a {\em singular} quadratic Lie algebra of {\em type
    $\Sb_3$} if $\dup(\g) = 3$.

\end{itemize}

\end{enumerate}

\end{defn}

Now, given a non-Abelian $n$-dimensional quadratic Lie algebra $\g$,
we can assume, up to i-isomorphism, that $\g = \CC^n$ equipped with
its canonical bilinear form $B$ (as a quadratic space) (\ref{1.3}). So
we introduce the following sets:

\begin{defn} For $n \geq 1$:

\begin{enumerate}

\item $\Qs(n)$ is the set of non-Abelian quadratic Lie algebra
  structures on $\CC^n$.

\item $\Os(n)$ is the set of {\em ordinary} quadratic Lie algebra
  structures on $\CC^n$.

\item $\Ss(n)$ is the set of {\em singular} quadratic Lie algebra
  structures on $\CC^n$.

\end{enumerate}

\end{defn}

By (\ref{1.5}), there is a one to one map from $\Qs(n)$ onto the
subset \[ \left\{ I \in \Wedge{}^3(\CC^n) \mid I \neq 0, \{I, I\} = 0
\right\} \subset \Wedge{}^3(\CC^n).\] In the sequel, we identify these
two sets, so that $\Qs(n) \subset \Wedge^3(\CC^n)$.

\begin{prop} \label{2.6}  One has:

\begin{enumerate}

\item $\Qs(n)$ is an affine variety in $\Wedge^3(\CC^n)$.

\item $\Os(n)$ is a Zariski-open subset of $\Qs(n)$.

\item $\Ss(n)$ is a Zariski-closed subset of $\Qs(n)$.

\end{enumerate}

\end{prop}

\begin{proof}
  The map $I \mapsto \{I,I\}$ is a polynomial map from
  $\Wedge^3(\CC^n)$ into $\Wedge^4(\CC^n)$, so the first claim
  follows.

  Fix $I \in \Wedge^3(\CC^n)$ such that $\{I, I\} = 0$. Consider the
  map ${\mathsf m} : (\CC^n)^* \to \Wedge^4(\CC^n)$ defined by
  ${\mathsf m}(\alpha) = \alpha \wedge I$, for all $\alpha \in
  (\CC^n)^*$. Then, if $\g$ is the quadratic Lie algebra associated to
  $I$, one has $\dup(\g) = 0$ if, and only if, $\rank({\mathsf m}) =
  n$. This can never happen for $n \leq 4$. Assume $n\geq 5$. Let $M$
  be a matrix of ${\mathsf m}$ and $\Delta_i$ be the minors of order
  $n$, for $1 \leq i\leq \binom{n}{4}$ . Then $\g \in \Os(n)$ if, and
  only if, there exists $i$ such that $\Delta_i \neq 0$. But
  $\Delta_i$ is a polynomial function and from that the second and the
  third claims follow.
\end{proof}

\begin{lem} \label{2.7} Let $\g_1$ and $\g_2$ be non-Abelian quadratic
  Lie algebras. Then $\g_1 \oplusp \g_2$ is an ordinary quadratic Lie
  algebra.
\end{lem}

\begin{proof}
  Set $\g = \g_1 \oplusp \g_2$. Denote by $I$, $I_1$ and $I_2$ the
  non-trivial 3-forms associated to $\g$, $\g_1$ and $\g_2$
  respectively. 

  One has $\Wedge(\g) = \Wedge(\g_1) \otimes \Wedge(\g_2)$,
  $\Wedge^k(\g) = \oplus_{r+s=k} \Wedge^r(\g_1) \otimes
  \Wedge^s(\g_2)$ and $I = I_1 + I_2$, with $I_1 \in \Wedge^3(\g_1)$
  and $I_2 \in \Wedge^3(\g_2)$. It immediately results that for
  $\alpha = \alpha_1 + \alpha_2 \in \g_1^* \oplus \g_2^*$, one has
  $\alpha \wedge I =0$ if, and only if, $\alpha_1 = \alpha_2 = 0$.
\end{proof}

\begin{prop} \label{2.6a} One has:

\begin{enumerate}

\item $\Qs(n) \neq \emptyset$ if, and only if, $n \geq 3$.

\item $\Os_3 = \Os_4 = \emptyset$ and $\Os(n) \neq \emptyset$ if $n
  \geq 6$.

\end{enumerate}

\end{prop}

\begin{proof}
  If $\g$ is a non-Abelian quadratic Lie algebra, using Remark
  \ref{2.3}, one has $\dim([\g, \g]) \geq 3$, so $\Qs(n) = \emptyset$
  if $n <3$.

  We shall now use some elementary quadratic Lie algebras given in
  Section 6 of \cite{PU07}. We denote these algebras by $\g_i$,
  according to their dimension, so that $\dim(\g_i) = i$, for $3 \leq
  i \leq 6$. Note that $\g_3 = \ok(3)$, $\g_4$, $\g_5$ and $\g_6$ are
  examples of elements of $\Qs(3)$, $\Qs(4)$, $\Qs(5)$ and $\Qs(6)$,
  respectively.

  Consider \[\g := \underset{3 \leq i \leq
    6}{\overset{\perp}{\bigoplus}} (\overbrace{\g_i \oplusp \dots
    \oplusp \g_i}^{k_i \textrm{times}}).\] Then $\dim(\g) =
  \sum_{i=3}^6 i k_i$ and by Lemma \ref{2.7}, $\dup(\g) = 0$, so we
  obtain $\Os(n) \neq \emptyset$ if $n \geq 6$.

  Finally, let $\g$ be a non-Abelian quadratic Lie algebra of
  dimension 3 or 4 with associated 3-form $I$. Then $I$ is
  decomposable, so $\g$ is singular. Therefore $\Os_3$ = $\Os_4 =
  \emptyset$.

\end{proof}
 
\begin{rem}
  We shall prove in Appendix 2 that $\Os_5 = \emptyset$. So,
  generically a non-Abelian quadratic Lie algebra is ordinary if $n
  \geq 6$.
\end{rem}

\begin{defn}
  A quadratic Lie algebra $\g$ is {\em indecomposable} if $\g = \g_1
  \oplusp \g_2$, with $\g_1$ and $\g_2$ ideals of $\g$, imply $\g_1$
  or $\g_2 = \{0\}$.
\end{defn}

The Proposition below gives another characterization of reduced
singular quadratic Lie algebras.

\begin{prop}\label{2.6b}
  Let $\g$ be a singular quadratic Lie algebra. Then $\g$ is reduced
  if, and only if, $\g$ is indecomposable.
\end{prop}

\begin{proof}
  If $\g$ is indecomposable, by Proposition \ref{2.8}, $\g$ is
  reduced. If $\g$ is reduced and $\g = \g_1 \oplusp \g_2$, with
  $\g_1$ and $\g_2$ ideals of $\g$, then $\Zs(\g_i) \subset [\g_i,
  \g_i]$ for $i =1,2$. So $\g_i$ is reduced or $\g_i = \{0\}$. But if
  $\g_1$ and $\g_2$ are both reduced, by Lemma \ref{2.7}, one has
  $\dup(\g) = 0$. Hence $\g_1$ or $\g_2 = \{0\}$.
\end{proof}

\section{Quadratic Lie algebras of type $\Sb_1$} \label{Section3}

\subsection{} \label{3.1} Let $(\g,B)$ be a quadratic vector space and $I$ be a
non-zero 3-form in $\Wedge^3(\g)$. As in (\ref{1.5}), we define a Lie
bracket on $\g$ by:
\[ [X,Y] := \phi^{-1} (\io_{X \wedge Y} (I)), \ \forall \ X, Y \in
\g.\] Then $\g$ becomes a quadratic Lie algebra with an invariant
bilinear form $B$ if, and only if, $\{I, I \} = 0$ \cite{PU07}.

In the sequel, we assume that $\dim(\Vs_I) = 1$. Fix $\alpha \in
\Vs_I$ and choose $\Omega \in \Wedge^2(\g)$ such that $I = \alpha
\wedge \Omega$ as follows: let $\{ \alpha, \alpha_1, \dots,
\alpha_r\}$ be a basis of $\Ws_I$. Then, $I \in \Wedge^3(\Ws_I)$ by
Proposition \ref{1.1}. We set:\[ X_0 := \phi^{-1}(\alpha) \ \text{ and
} X_i := \phi^{-1}(\alpha_i), \ 1\leq i \leq r.\] So, we can choose
$\Omega \in \Wedge^2(V)$ where $V = \spa\{X_1, \dots, X_r \}$.  Note
that $\Omega$ is an indecomposable bilinear form, so $\dim(V) > 3$.

We define $C : \g \to \g$ by \[ B(C(X), Y) := \Omega(X,Y).\] Therefore
$C$ is skew-symmetric with respect to $B$.

\begin{lem} \label{3.2} The following are equivalent:

\begin{enumerate}

\item $\{I,I\} = 0$ 

\item $\{\alpha, \alpha\} = 0$ and $\{\alpha, \Omega\} = 0$

\item $B(X_0, X_0) =0$ and $C(X_0) = 0$

\end{enumerate}

In this case, one has $\dim([\g,\g]) > 4$, $\Zs(\g) \subset \ker(C)$,
$\im(C) \subset [\g,\g]$ and $X_0 \in \Zs(\g) \cap [\g, \g]$.
\end{lem}

\begin{proof}
  It is easy to see that: \[ \{I, I\} = 0 \Leftrightarrow \{\alpha,
  \alpha \} \wedge \Omega \wedge \Omega = 2 I \wedge \{\alpha, \Omega
  \}.\] If $\Omega \wedge \Omega = 0$, then $\Omega$ is decomposable
  and that is a contradiction since $\dim(\Vs_I) = 1$. So $ \Omega
  \wedge \Omega \neq 0$.

  If $\{ \alpha, \alpha \} \neq 0$, then $\alpha$ divides $\Omega
  \wedge \Omega \in \Wedge^4(V)$, another contradiction. That implies
  $\{\alpha, \alpha \} = 0 = B(X_0, X_0)$. It results that $\{ \alpha,
  \Omega \} \in \Vs_I = \CC \alpha$, hence $\{\alpha, \Omega \} =
  \lambda \alpha$ for some $\lambda \in \CC$. But $\{ \alpha, \Omega
  \}$ is an element of $\Wedge^1(V)$, so $\lambda$ must be zero and by
  (\ref{1.4}), $\io_{X_0} (\Omega) = 0$, therefore $C(X_0) =
  0$. Moreover, since $\{ \alpha, \alpha \} = \{ \alpha, \Omega \} =
  0$, using $I = \alpha \wedge \Omega$, we deduce that $\{\alpha, I \}
  = 0$. Again by (\ref{1.4}), it results that $B(X_0, [X,Y]) = \{
  \alpha, I \} (X \wedge Y) = 0$, for all $X$, $Y \in \g$. So $X_0 \in
  [\g,\g]^\perp = \Zs(\g)$. Also, $\Vs_I \subset \Ws_I$, so $X_0 =
  \phi^{-1} (\alpha) \in \phi^{-1} (\Ws_I) = [\g,\g]$.

  Write $\Omega = \sum_{i<j} a_{ij} \alpha_i \wedge \alpha_j$, with
  $a_{ij} \in \CC$. Since $\Ws_I = \phi([\g, \g])$ and $X_1, \dots,
  X_r \in [\g, \g]$, we deduce that \[ C = \sum_{i<j} a_{ij} (\alpha_i
  \otimes X_j - \alpha_j\otimes X_i)\] Hence $\im(C) \subset
  [\g,\g]$. Since $C$ is skew-symmetric, one has $\ker(C) =
  \im(C)^\perp$ and it follows $\Zs(\g) = [\g, \g]^\perp \subset
  \ker(C)$.


  Finally, $[\g, \g] = \CC X_0 \oplus V$ and since $\dim(V) > 3$, we
  conclude that $\dim([\g,\g]) > 4$.
\end{proof}

\begin{rem} \label{3.3}
  It is important to notice that our choice of $\Omega$ such that $I =
  \alpha \wedge \Omega$ is not unique, it depends on the choice of
  $V$, so $C$ is not uniquely defined. Assume we consider another
  vector space $V'$ and $I = \alpha \wedge \Omega'$. Then $\Omega' =
  \Omega + \alpha \wedge \beta$ for some $\beta \in \g^*$. Let $X_1 =
  \phi^{-1}(\beta)$ and let $C'$ be the map associated to
  $\Omega'$. By a straightforward computation, $C' = C + \alpha
  \otimes X_1 - \beta \otimes X_0$. Since $C'(X_0) = 0$, we must have
  $B(X_0, X_1) = 0$.
\end{rem}

\subsection{} We keep the notation as in the previous
subsection. Assume that $\{I, I \}= 0$. Hence $\g$ is a quadratic Lie
algebra of type $\Sb_1$.

\begin{lem} \label{3.4}
  There exists $Y_0 \in V^\perp$ such that \[ V^\perp = \Zs(\g) \oplus
  \CC Y_0,\ B(Y_0, Y_0) = 0 \ \text{ and } \ B(X_0, Y_0) = 1.\]
  Moreover \[C(Y_0) = 0.\]
\end{lem}

\begin{proof}
  One has $\phi^{-1} (\Ws_I) = [\g,\g] = \CC X_0 \oplus V$, therefore
  $\Zs(\g) \subset V^\perp$ and $\dim(\Zs(\g)) = \dim(\g) -
  \dim([\g,\g]) = \dim(V^\perp) - 1$. So there exists $Y \in V^\perp$
  such that $V^\perp = \Zs(\g) \oplus \CC Y$. Now, $Y$ cannot be
  orthogonal to $X_0$, since it would be orthogonal to $[\g,\g]$ and
  therefore an element of $\Zs(\g)$. So we can assume that $B(X_0, Y)
  =1$. Replace $Y$ by $Y_0 = Y - \dfrac12 B(Y,Y) X_0$ to obtain
  $B(Y_0, Y_0) = 0$ (recall $B(X_0,X_0) = 0$).

  By Lemma \ref{3.2}, $\im(C) \subset V$ and that implies $B(Y_0,
  C(X)) = - B(C(Y_0),X) = 0$, for all $X \in \g$. Then $C(Y_0) = 0$.
\end{proof}

\begin{prop} \label{3.5}
We keep the previous notation and assumptions. Then:

\begin{enumerate}

\item $[X,Y] = B(X_0, X) C(Y) - B(X_0,Y) C(X) + B(C(X),Y) X_0$, for
  all $X$, $Y \in \g$.

\smallskip

\item $C = \ad(Y_0)$ and $\rank(C)$ is even.

\smallskip

\item $\ker(C) = \Zs(\g) \oplus \CC Y_0$, $\im(C) = V$ and $[\g,\g] =
  \CC X_0 \oplus \im(C)$.

\smallskip

\item the Lie algebra $\g$ is solvable. Moreover, $\g$ is nilpotent
  if, and only if, $C$ is nilpotent.

\smallskip

\item the dimension of $[\g,\g]$ is greater or equal to 5 and it is
  odd. 
\end{enumerate}

\end{prop}

\begin{proof} \hfill

\begin{enumerate}

\item This is a straightforward computation, use $B([X,Y],Z) = (\alpha
  \wedge \Omega) (X,Y,Z)$, $\alpha(X) = B(X_0,X)$ and $\Omega(X,Y) =
  B(C(X),Y)$, for all $X, Y, Z \in \g$.

\item Set $X = Y_0$ in (1) and use Lemma \ref{3.4} to show $C =
  \ad(Y_0)$. Since $C(\g) = \ad(Y_0)(\g) = \phi^{-1} \left(
    \check{\ad}(\g) (\phi(Y_0) )\right)$, the rank of $C$ is the
  dimension of the coadjoint orbit through $\phi(Y_0)$, so it is even
  (see also Appendix 1).

\item We may assume that $\g$ is reduced. Then $\Zs(\g)$ is totally
  isotropic and $\Zs(\g) \subset X_0^\perp$.  Write $X_0^\perp =
  \Zs(\g) \oplus \hk$ with $\hk$ a complementary subspace of
  $\Zs(\g)$. Therefore $\g = \Zs(\g) \oplus \hk \oplus \CC Y_0$ and
  for an element $X = Z + H + \lambda Y_0 \in \ker(C)$, we deduce $H
  \in \ker(C)$ by Lemmas \ref{3.2} and \ref{3.4}.

  But $B(X_0,H) = 0$, so using (1), $H \in \Zs(\g)$. It results that
  $H = 0$. Then $\ker(C) = \Zs(\g) \oplus \CC Y_0$. In addition,
  \[ \dim(\im(C)) = \dim(\hk) = \dim(X_0^\perp) - \dim(\Zs(\g)) =
  \dim([\g,\g]) - 1.\]

  Our choice of $V$ implies that $[\g,\g] = \phi^{-1} (\Ws_I) = \CC
  X_0 \oplus V$ and $\im(C) \subset V$ (see the proof of Lemma
  \ref{3.2}). Therefore $\im(C) = V$ and $[\g,\g] = \CC X_0 \oplus
  \im(C)$.

\item Since $B(X_0, \im(C)) = 0$, then $[[\g,\g], [\g,\g]] = [\im(C),
  \im(C)] \subset \CC X_0$. We conclude that $\g$ is solvable. If $\g$
  is nilpotent, then $C = \ad(Y_0)$ is nilpotent. If $C$ is nilpotent,
  using $\im(C) \subset X_0^\perp$, we obtain by induction that
  $(\ad(X))^k(\g) \subset \CC X_0 \oplus \im(C^k)$ for any $k \in
  \NN$. So $\ad(X)$ is nilpotent, for all $X \in \g$ and that implies
  $\g$ nilpotent.

\item Notice that $[\g,\g] = \CC X_0 \oplus \im(C)$ and $\rank(C)$ is
  even, so $\dim([\g,\g])$ is odd. By (\ref{3.1}), $\dim([\g,\g]) \geq
  5$.
\end{enumerate}

\end{proof}

\subsection{} Recall that $C$ is not unique (see Remark \ref{3.3}) and
it depends on the choice of $V$. Let \[ \ak := X_0^\perp / \CC X_0.\]
We denote by $\hat{X}$ the class of an element $X \in \g$.

\begin{prop} \label{3.6} \hfill

Keep the notation above. One has:

\begin{enumerate}

\item the Lie algebra $\ak$ is Abelian.

\item Define \[\hat{B} (\hat{X}, \hat{Y}) := B(X,Y), \ \forall \ X, Y
  \in \g.\] Then $\hat{B}$ is a non degenerate symmetric bilinear form
  on $\ak$.

\smallskip

\item Define \[\hat{C} (\hat{X}) := C(X), \ \forall \ X \in \g.\] Then
  $\hat{C} \in \Lc(\ak)$ is a skew-symmetric map with $\rank(\hat{C})
  = \rank(C)$ even and $\rank(\hat{C}) \geq 4$.

\smallskip

\item $\hat{C}$ does not depend on the choice of $V$. More precisely,
  if $\Ws_I = \CC \alpha \oplus \phi(V')$ and $C'$ is the associated
  map to $V'$ (see Remark \ref{3.3}), then $\hat{C'}= \hat{C}$.

\smallskip

\item the Lie algebra $\g$ is reduced if, and only if, $\ker(\hat{C})
  \subset \im(\hat{C})$.

\end{enumerate}

\end{prop}

\begin{proof} \hfill

\begin{enumerate}

\item It follows from Proposition \ref{3.5} (1).

\item It is clear that $\hat{B}$ is well-defined. Now, since
  $B(X_0,Y_0) = 1$, $B(X_0,X_0) = B(Y_0,Y_0) = 0$, the restriction of
  $B$ to $\spa\{X_0, Y_0 \}$ is non degenerate. So \[\g = \spa\{X_0,
  Y_0 \} \oplusp \spa\{X_0, Y_0 \}^\perp,\] $X_0^\perp = \CC X_0
  \oplus \spa\{X_0, Y_0 \}^\perp$ and $X_0^\perp{}^\perp = X_0^\perp
  \cap \spa\{X_0, Y_0 \} = \CC X_0$. We conclude that $\hat{B}$ is non
  degenerate.

\item We have $C(X_0^\perp) = \ad(Y_0)(X_0^\perp) \subset X_0^\perp$
  since $X_0^\perp$ is an ideal of $\g$. Moreover, $C(X_0) = 0$, so
  $\hat{C}$ is well-defined. The image of $C$ is contained in
  $X_0^\perp$ and $\im(C) \cap \CC X_0 = \{0\}$, therefore $\dim
  \left(\im(C)/\CC X_0\right) = \dim(\im(\hat{C})) =
  \dim(\im(C))$. Now it is enough to apply Proposition \ref{3.5}.

\item By Remark \ref{3.3}, we have $C' = C + \alpha \otimes X_1 -
  \beta \otimes X_0$. But $\alpha(X_0) = 0$, so $\hat{C'} = \hat{C}$.

\item By Proposition \ref{3.5}, we have $\ker(C) = \Zs(\g) \oplus \CC
  Y_0$ and by Lemma \ref{3.2}, we have $\Zs(\g) \subset
  X_0^\perp$. Again by Proposition \ref{3.5}, we conclude that
  $\ker(\hat{C} ) = \Zs(\g) / \CC X_0$. Applying Proposition \ref{3.5}
  once more, we have $[\g,\g] = \CC X_0 \oplus \im(C)$, so
  $\im(\hat{C}) = [\g, \g]/ \CC X_0$. Then $\ker(\hat{C}) \subset
  \im(\hat{C})$ if, and only if, $\Zs(\g) \subset [\g,\g] + \CC
  X_0$. But $X_0 \in [\g,\g]$ (see Lemma \ref{3.2}), so the result
  follows.

\end{enumerate}

\end{proof}

We should notice that $\hat{C}$ still depends on the choice of
$\alpha$ (see Remark \ref{3.3}): if we replace $\alpha$ by $\lambda
\alpha$, for a non-zero $\lambda \in \CC$, that will change $\hat{C}$
into $\dfrac{1}{\lambda} \hat{C}$. So there is not a {\em unique} map
$\hat{C}$ associated to $\g$ but rather a {\em family} $\{ \lambda
\hat{C} \mid \lambda \in \CC \setminus \{0\} \}$ of associated
maps. In other words, there is a line \[[\hat{C}] := \{ \lambda
\hat{C} \mid \lambda \in \CC \} \in \ps(\ok(\ak)) \] where
$\ps(\ok(\ak))$ is the projective space associated to the space
$\ok(\ak)$.

\begin{defn}
  We call $[\hat{C}]$ the {\em line of skew-symmetric maps} associated
  to the quadratic Lie algebra $\g$ of type $\Sb_1$.
\end{defn}

\begin{rem} \label{3.7} The unicity of $[\hat{C}]$ is valuable, but
  the fact that $\hat{C}$ acts on a quotient space and not on a
  subspace of $\g$ could be a problem. Hence it is convenient to use
  the following decomposition of $\g$: the restriction of $B$ to $\CC
  X_0 \oplus \CC Y_0$ is non degenerate, so we can write $\g = (\CC
  X_0 \oplus \CC Y_0) \oplusp \qk$ where $\qk = (\CC X_0 \oplus \CC
  Y_0)^\perp$. Since $C(X_0) = C(Y_0) = 0$ and $C \in \ok(\g)$, $C$
  maps $\qk$ into $\qk$. Let $\pi : X_0^\perp \to X_0^\perp / \CC X_0$
  be the canonical surjection and $\cb = C|_\qk$. Then the restriction
  $\pi_\qk : \qk \to X_0^\perp / \CC X_0$ is an isometry and $\hat{C}
  = \pi_\qk \ \cb \ \pi_\qk^{-1}$.

  Remark that $Y_0$ is not unique, but if $Y_0'$ satisfies Lemma
  \ref{3.4}, consider $C' = \ad(Y_0')$ and $\qk'$ such that $\g = (\CC
  X_0 \oplus \CC Y_0') \oplusp \qk'$, therefore $\hat{C} = \pi_\qk' \
  \cpb \ \pi_\qk'^{-1}$ with the obvious notation. It results that $
  \pi_\qk'^{-1} \pi_\qk$ is an isometry from $\qk$ to $\qk'$ and
  that \[\cpb = \left( \pi_\qk'^{-1} \pi_\qk \right) \cb \left(
    \pi_\qk'^{-1} \pi_\qk \right)^{-1}.\]

  We shall develop this aspect in the next Section.
\end{rem}

\section{Solvable singular quadratic Lie algebras and double
  extensions}

\subsection{} Double extensions are a very effective method initiated
by V. Kac to construct quadratic Lie algebras (see \cite{Kac, FS,
  MR}). Here, we only need a particular case that we shall recall:

\begin{defn} \hfill

\begin{enumerate}

\item Let $(\qk, B_\qk)$ be a quadratic vector space and $\cb : \qk
  \to \qk$ be a skew-symmetric map. Let $(\tk = \spa \{ X_1, Y_1 \},
  B_\tk)$ be a 2-dimensional quadratic vector space with $B_\tk$
  defined by \[ B_\tk(X_1, X_1) = B_\tk(Y_1, Y_1) = 0, \ B_\tk(X_1,
  Y_1) = 1.\] Consider \[\g = \qk \oplusp \tk\] equipped with a
  bilinear form $B := B_\qk + B_\tk$ and define a bracket on $\g$
  by \[ [X + \lambda X_1 + \mu Y_1, Y + \lambda' X_1 + \mu' Y_1] :=
  \mu \cb (Y) - \mu' \cb (X) + B( \cb (X), Y) X_1, \] for all $X, Y
  \in \qk, \lambda, \mu, \lambda', \mu' \in \CC$. Then $(\g, B)$ is a
  quadratic solvable Lie algebra. We say that $\g$ is the {\em double
    extension} of $\qk$ by $\cb$.

\item Let $\g_i$ be double extensions of quadratic vector spaces
  $(\qk_i, B_i)$ by skew-symmetric maps $\cb_i \in \Lc(\qk_i)$, for $1
  \leq i \leq k$. The {\em amalgamated product} \[ \g = \g_1 \ta \g_2
  \ta \dots \ta \g_k\] is defined as follows:

\begin{itemize}

\item consider $(\qk,B)$ be the quadratic vector space with $\qk =
  \qk_1 \oplus \qk_2 \oplus \dots \oplus \qk_k$ and the bilinear form
  $B$ such that $B(\sum_{i = I}^k X_i, \sum_{i = I}^k Y_i) = \sum_{i =
    I}^k B_i(X_i, Y_i)$, for $X_i, Y_i \in \qk_i$, $1 \leq i \leq k$.

\item the skew-symmetric map $\cb \in \Lc(\qk)$ is defined by $\cb(
  \sum_{i = I}^k X_i) =$ \linebreak $ \sum_{i = I}^k \cb_i(X_i)$, for
  $X_i \in \qk_i$, $1 \leq i \leq k$.

\end{itemize}
 
\medskip

Then $\g$ is the double extension of $\qk$ by $\cb$.

\end{enumerate}

\end{defn}

In this Section, we will show that double extensions are highly
related to singular quadratic Lie algebras. Amalgamated products will
be used in Sections 5 and 6 to {\bf decompose} double extensions.

We notice that if $\g_1 \iiso \g_1'$ and $\g_2 \iiso \g_2'$, it may
happen that $\g_1 \ta \g_2$ and $\g_1' \ta \g_2'$ are not even
isomorphic. So, amalgamated products have a bad behavior with respect
to i-isomorphisms. An example will be given in Section 5, Remark
\ref{5.12}.

\begin{lem} \label{4.3} We keep the notation above. 

\begin{enumerate}

\item Let $\g$ be the double extension of $\qk$ by
  $\cb$. Then \[ [X, Y] = B(X_1, X) C (Y) - B(X_1, Y) C(X) +
  B(C(X), Y) X_1, \ \forall \ X, Y \in \g, \] where $C =
  \ad(Y_1)$. Moreover, $X_1 \in \Zs(\g)$ and $C|_\qk = \cb$.

\item Let $\g'$ be the double extension of $\qk$ by $\cpb =
  \lambda \cb$, $\lambda \in \CC$, $\lambda \neq 0$. Then
  $\g$ and $\g'$ are i-isomorphic.

\end{enumerate}

\end{lem}

\begin{proof}\hfill

\begin{enumerate}

\item This is a straightforward computation.

\item Write $\g = \qk \oplusp \tk = \g'$. Denote by $[\cdot, \cdot]'$
  the Lie bracket on $\g'$. Define $A : \g \to \g'$ by $A(X_1) =
  \lambda X_1$, $A(Y_1) = \dfrac{1}{\lambda} Y_1$ and $A|_\qk =
  \Id_\qk$. Then $A([Y_1,X]) = C(X) = [A(Y_1), A(X)]'$ and $A([X,Y]) =
  [A(X), A(Y)]'$, for all $X, Y \in \qk$. So $A$ is an i-isomorphism.

\end{enumerate}

\end{proof}

\subsection{} A natural consequence of formulas in Lemma \ref{4.3} and
Proposition \ref{3.5} (1) is given by the Proposition below:

\begin{prop} \label{4.4} \hfill

\begin{enumerate}

\item Consider the notation in Section \ref{Section3}, Remark
  \ref{3.7}. Let $\g$ be quadratic Lie algebra of type $\Sb_1$ (that
  is, $\dup(\g) = 1$). Then $\g$ is the double extension of $\qk =
  (\CC X_0 \oplus \CC Y_0)^\perp$ by $\cb = \ad(Y_0)|_\qk$.

\smallskip

\item Let $(\g,B)$ be a quadratic Lie algebra. Let $\g'$ be a double
  extension of a quadratic vector space $(\qk',B')$ by a map
  $\cpb$. Let $A$ be an i-isomorphism of $\g'$ onto $\g$ and
  write $\qk = A(\qk')$. Then $\g$ is a double extension of $(\qk,
  B|_{\qk \times \qk})$ by the map $\cb = \overline{A} \
  \cpb \ \overline{A}^{-1}$ where $\overline{A} = A|_{\qk'}$.

\smallskip

\item Let $\g$ be the double extension of a quadratic vector space
  $\qk$ by a map $\cb \neq 0$. Then $\g$ is a singular
  solvable quadratic Lie algebra. Moreover:

\begin{itemize}

\item[(a)] $\g$ is of type $\Sb_3$ if, and only if,
  $\rank(\cb) = 2$.

\smallskip

\item[(b)] $\g$ is of type $\Sb_1$ if, and only if,
  $\rank(\cb) \geq 4$.

\smallskip

\item[(c)] $\g$ is reduced if, and only if, $\ker(\cb)
  \subset \im(\cb)$.

\smallskip

\item[(d)] $\g$ is nilpotent if, and only if, $\cb$ is nilpotent.

\end{itemize}

\end{enumerate}

\end{prop}

\begin{proof} \hfill

\begin{enumerate}

\item Let $\bk = \CC X_0 \oplus \CC Y_0$. Then $B|_{\bk \times \bk}$
  is non degenerate and $\g = \bk \oplus \qk$. Since $\ad(Y_0)(\bk)
  \subset \bk$ and $\ad(Y_0)$ is skew-symmetric, we have
  $\ad(Y_0)(\qk) \subset \qk$. By Proposition \ref{3.5} (1), we
  have \[ [X, X'] = B ( \cb(X), X') X_0, \ \forall \ X, X' \in \qk.\]
  Set $X_1 := X_0$ and $Y_1 := Y_0$ to obtain the result.

\item Write $\g' = (\CC X_1' \oplus \CC Y_1') \oplusp \qk'$. Let $X_1
  = A(X_1')$ and $Y_1 = A (Y_1')$. Then $\g = (\CC X_1 \oplus \CC Y_1)
  \oplusp \qk$ and \[ [Y_1, X] = (A \cpb A^{-1}) (X), \
  \forall \ X \in \qk, \text{ and } \] \[ [X, Y] = B((A \cpb
  A^{-1}) (X), Y) X_1', \ \forall \ X, Y \in \qk.\] and this proves
  the result.

\item Let $\g = (\CC X_1 \oplus \CC Y_1) \oplusp \qk$, $C = \ad(Y_1)$,
  $\alpha = \phi(X_1)$, $\Omega(X,Y) = B(C(X),Y)$, for all $X$, $Y \in
  \g$ and $I$ be the 3-form associated to $\g$. Then the formula for
  the Lie bracket in Lemma \ref{4.3}(1) can be translated as $I =
  \alpha \wedge \Omega$, hence $\dup(\g) \geq 1$ and $\g$ is singular.

  Let $W_\Omega$ be the set $W_\Omega = \{ \io_X(\Omega), X \in \g
  \}$. Then $W_\Omega = \phi(\im(\cb))$. Therefore
  $\rank(\cb) \geq 2$ by Proposition \ref{1.1} and $\Omega$
  is decomposable if, and only if, $\rank(\cb) =2$.

  If $\rank(\cb) > 2$, then $\g$ is of type $\Sb_1$ and by
  Proposition \ref{3.6}, we have $\rank(\cb) \geq 4$.

  Finally, $\Zs(\g) = \CC X_1 \oplus \ker(\cb)$ and $[\g,\g]
  = \CC X_1 \oplus \im(\cb)$, so $\g$ is reduced if, and only
  if, $\ker(\cb) \subset \im(\cb)$.

  The proof of the last claim is exactly the same as in Proposition
  \ref{3.5} (4).

\end{enumerate}

\end{proof}

\subsection{} A complete classification (up to i-isomorphism) of
quadratic Lie algebras of type $\Sb_3$ is given in \cite{PU07}. We
shall recall the characterization of these algebras here:

\begin{prop} \label{4.5} Let $\g$ be a quadratic Lie algebra of type
  $\Sb_3$. Then $\g$ is i-isomorphic to an algebra $ \lk \oplusp \zk$
  where $\zk$ is a central ideal of $\g$ and $\lk$ is one of the
  following algebras:

\begin{enumerate}

\item $\g_3(\lambda) = \ok(3)$ equipped with the bilinear form $B =
  \lambda \kappa$ where $\kappa$ is the Killing form and $\lambda \in
  \CC$, $\lambda \neq 0$.

\smallskip

\item $\g_4$, a 4-dimensional Lie algebra: consider $\qk = \CC^2$,
  $\{E_1, E_2 \}$ its canonical basis and the bilinear form $B$
  defined by $B(E_1,E_1) = B(E_2, E_2) = 0$ and $B(E_1, E_2) =
  1$. Then $\g_4$ is the double extension of $\qk$ by the
  skew-symmetric map \[ \cb = \begin{pmatrix} 1 & 0 \\ 0 &
    -1 \end{pmatrix}.\] Moreover, $\g_4$ is solvable, but it is not
  nilpotent.

\smallskip

\item $\g_5$, a 5-dimensional Lie algebra: consider $\qk = \CC^3$,
  $\{E_1, E_2, E_3 \}$ its canonical basis and the bilinear form $B$
  defined by $B(E_1,E_1) = B(E_2, E_2) = B(E_1, E_2)$ $ = B(E_2, E_3)
  = 0$ and $B(E_1, E_3) = B(E_2, E_2) = 1$. Then $\g_5$ is the double
  extension of $\qk$
  by the skew-symmetric map \[ \cb = \begin{pmatrix} 0 & 1 & 0 \\
    0 & 0 & -1 \\ 0 & 0 & 0 \end{pmatrix}.\] Moreover, $\g_5$ is
  nilpotent.

\smallskip

\item $\g_6$, a 6-dimensional Lie algebra: consider $\qk = \CC^4$,
  $\{E_1, E_2, E_3, E_4 \}$ its canonical basis and the bilinear form
  $B$ defined by $B(E_1,E_3) = B(E_2, E_4) = 1$ and $B(E_i, E_j) = 0$
  otherwise. Then $\g_6$ is the double extension of $\qk$
  by the skew-symmetric map \[ \cb = \begin{pmatrix} 0 & 1 & 0 & 0\\
    0 & 0 & 0 & 0 \\ 0 & 0 & 0 & 0 \\ 0 & 0 & -1 & 0 \end{pmatrix}.\]
  Moreover, $\g_6$ is nilpotent.

\end{enumerate}

\smallskip

{\bf {\em All solvable quadratic Lie algebras of type $\Sb_3$ are
    double extensions of a quadratic vector space by a skew-symmetric
    map.}}

\end{prop}

\subsection{} Let $(\qk,B)$ be a quadratic vector space. We recall
that $\OO(\qk)$ is the group of orthogonal maps and $\ok(\qk)$ is its
Lie algebra, i.e. the Lie algebra of skew-symmetric maps. Recall that
the {\em adjoint action} is the action of $\OO(\qk)$ on $\ok(\qk)$ by
conjugation.

\begin{prop} \label{4.6} Let $(\qk, B)$ be a quadratic vector
  space. Let $\g = (\CC X_1 \oplus \CC Y_1) \oplusp \qk$ and $\g' =
  (\CC X_1' \oplus \CC Y_1') \oplusp \qk$ be double extensions of
  $\qk$, by skew-symmetric maps $\cb$ and $\cpb$ respectively. Then:

\begin{enumerate}

\item there exists a Lie algebra isomorphism between $\g$ and $\g'$
  if, and only if, there exists an invertible map $P \in \Lc(\qk)$ and
  a non-zero $\lambda \in \CC$ such that $\cpb = \lambda \ P
  \cb P^{-1}$ and $P^* P \cb = \cb$, where
  $P^*$ is the adjoint map of $P$ with respect to $B$.

\smallskip

\item there exists an i-isomorphism between $\g$ and $\g'$ if, and
  only if, $\cpb$ is in the $\OO(\qk)$-adjoint orbit through
  $\lambda \cb$ for some non-zero $\lambda \in \CC$.

\end{enumerate}

\end{prop}

\begin{proof}  \hfill

\begin{enumerate}

\item Let $A : \g \to \g'$ be a Lie algebra isomorphism. We know by
  Proposition \ref{4.4} that $\g$ and $\g'$ are singular. Assume that
  $\g$ is of type $\Sb_3$. Then $3 = \dim( [\g,\g]) = \dim(
  [\g',\g'])$. So $\g'$ is also of type $\Sb_3$
  (\cite{PU07}). Therefore, $\g$ and $\g'$ are either both of type
  $\Sb_1$ or both of type $\Sb_3$. Let us study these two cases.

\begin{itemize}

\item[(i)] First, assume that $\g$ and $\g'$ are both of type
  $\Sb_1$. We start by proving that $A (\CC X_1 \oplus \qk) = \CC X_1'
  \oplus \qk$. If this is not the case, there is $X \in \qk$ such that
  $A(X) = \beta X_1' + \gamma Y_1' + Y$ with $Y \in \qk$ and $\gamma
  \neq 0$. Then \[ [ A(X), \CC X_1' \oplus \qk]' = \gamma \cpb(\qk) +
  [Y,\qk]'.\] Since $\g'$ is of type $\Sb_1$, we have $\rank(\cb')
  \geq 4$ (see Proposition \ref{4.4}) and it follows that $\dim \left(
    [ A(X), \CC X_1' \oplus \qk]' \right) \geq 4$. On the other hand,
  $ [ A(X), \CC X_1' \oplus \qk]'$ is contained in $A([X,\g])$ and
  $\dim( [X,\g]) \leq 2$, so we obtain a contradiction.

  Next, we prove that $A(X_1) \in \CC X_1'$. Since $X_1 \in [\g,\g]$,
  then there exists $X$, $Y \in \qk$ such that $X_1 = [X,Y]$. Then $A
  (X_1) = [A(X), A(Y)]' \in [ \CC X_1' \oplus \qk, \CC X_1' \oplus
  \qk]' = \CC X_1'$. Hence $A(X_1) = \mu X_1'$ for some non-zero $\mu
  \in \CC$.

  Now, write $A|_\qk = P + \beta \otimes X_1'$ with $P : \qk \to \qk$
  and $\beta \in \qk^*$. If $X\in \ker(P)$, then$ A \left( X -
    \dfrac{1}{\mu} \beta(X) X_1 \right) = 0,$ so $X = 0$ and
  therefore, $P$ is invertible.

  For all $X$, $Y \in \qk$, we have $A ([X,Y]) = \mu B(\cb(X), Y)
  X_1'$. Also,
\begin{eqnarray*}
  A([X,Y]) &=& [P(X) + \beta(X)
  X_1', P(Y) + \beta(Y) X_1']' \\ &=& B ( \cpb P(X), P(Y)) X_1'.
\end{eqnarray*}
So it results that $P^* \cpb P = \mu \cb$.

Moreover, $A([Y_1, X]) = P(C(X) + \beta(C(X)) X_1'$, for all $X \in
\qk$. Let $A(Y_1) = \gamma Y_1' + Y + \delta X_1'$, with $Y \in
\qk$. Therefore \[ A([Y_1, X]) = \gamma \cpb P(X) + B(\cpb (Y), P(X))
X_1'\] and we conclude that $P \ \cb \ P^{-1} = \gamma \cpb$ and since
$P^* \cpb P = \mu \cb$, then $P^* P \cb = \gamma \mu \cb$.

Set $Q = \dfrac{1}{ (\mu \gamma)^{\frac12}} P$. It follows that $Q \cb
Q^{-1} = \gamma \cpb$ and $Q^* Q \cb = \cb$. This finishes the proof
in the case $\g$ and $\g'$ of type $\Sb_1$.

\item[(ii)] We proceed to the case when $\g$ and $\g'$ of type
  $\Sb_3$: the proof is a straightforward case-by-case
  verification. By Proposition \ref{2.8}, we can assume that $\g$ and
  $\g'$ are reduced. Then $\dim(\qk) = 2, 3$ or $4$ by Proposition
  \ref{4.5}.

  Recall that $\g$ is nilpotent if, and only if, $\cb$ is
  nilpotent (see Proposition \ref{4.4} (3)). The same is valid for
  $\g'$.

  If $\dim(\qk) = 2$, then $\g$ is not nilpotent, so $\cb$ is not
  nilpotent, $\Tr(\cb) = 0$ and $\cb$ must be semi-simple. Therefore
  we can find a basis $\{e_1, e_2 \}$ of $\qk$ such that $B(e_1, e_2)
  = 1$, $B(e_1, e_1) = B(e_2, e_2) = 0$ and the matrix of $\cb$ is
  $\begin{pmatrix} \mu & 0 \\ 0 & - \mu \end{pmatrix}$. The same holds
  for $\cpb$: there exists a basis $\{e_1', e_2' \}$ of $\qk$ such
  that $B(e_1', e_2') = 1$ and $B(e_1', e_1)' = B(e_2', e_2') = 0$
  such that the matrix of $\cpb$ is $\begin{pmatrix} \mu' & 0 \\ 0 & -
    \mu' \end{pmatrix}$. It results that $\cpb$ and $\dfrac{\mu'}{\mu}
  \cb$ are $\OO(\qk)$-conjugate and we are done.

  If $\dim(\qk) = 3$ or 4, then $\g$ and $\g'$ are nilpotent. We use
  the classification of nilpotent orbits given for instance in
  \cite{CM}: there is only one orbit in dimension 3 or 4, so $\cb$ and
  $\cpb$ are conjugate by $\OO(\qk)$.

\end{itemize}

\smallskip

This finishes the proof of the necessary condition. To prove the
sufficiency, we replace $\cpb$ by $\lambda P \cb P^{-1}$ to obtain $P^*
\cpb P = \lambda \cb$. Then we define $A : \g \to \g'$ by $A(X_1) =
\lambda X_1'$, $A(Y_1) = \dfrac{1}{\lambda} Y_1'$ and $A(X) = P(X)$,
for all $X \in \qk$. By a direct computation, we have for all $X$ and
$Y \in \qk$: \[ A([X,Y]) = [A(X), A(Y)]' \ \text{ and } \ A([Y_1, X])
= [A(Y_1), A(X)]',\] so $A$ is a Lie algebra isomorphism between $\g$
and $\g'$.

\smallskip

\item If $\g$ and $\g'$ are i-isomorphic, then the isomorphism $A$ in
  the proof of (1) is an isometry. Hence $P \in \OO(\qk)$ and $P^*
  \cpb P = \mu \cb$ gives the result.

  Conversely, define $A$ as above (sufficiency of (1)). Then $A$ is an
  isometry and it is easy to check that $A$ is an i-isomorphism.
\end{enumerate}

\end{proof}

\begin{cor} \label{4.7} Let $(\g,B)$ and $(\g',B')$ be double
  extensions of $(\qk, \overline{B})$ and $(\qk', \overline{B'})$
  respectively, where $\overline{B} = B|_{\qk \times \qk}$ and
  $\overline{B'} = B'|_{\qk' \times \qk'}$. Write $\g = (\CC X_1
  \oplus \CC Y_1) \oplusp \qk$ and $\g' = (\CC X_1' \oplus \CC Y_1')
  \oplusp \qk'$. Then:

\begin{enumerate}

\item there exists an i-isomorphism between $\g$ and $\g'$ if, and
  only if, there exists an isometry $\overline{A} : \qk \to \qk'$ such
  that $\cpb = \lambda \ \overline{A} \ \cb \ \overline{A}^{-1}$, for
  some non-zero $\lambda \in \CC$.

\smallskip

\item there exists a Lie algebra isomorphism between $\g$ and $\g'$
  if, and only if, there exist invertible maps $\overline{Q} : \qk \to
  \qk'$ and $\overline{P} \in \Lc(\qk)$ such that 

\begin{itemize}

\item[(i)] $\cpb = \lambda \ \overline{Q} \ \cb \ \overline{Q}^{-1}$
  for some non-zero $\lambda \in \CC$, 

\item[(ii)] $\overline{P}^* \ \overline{P} \ \cb = \cb$ and

\item[(iii)] $\overline{Q} \ \overline{P}^{-1}$ is an isometry from
  $\qk$ onto $\qk'$.

\end{itemize}

\end{enumerate}

\end{cor}

\begin{proof} \hfill

\begin{enumerate}

\item We can assume that $\dim(\g) = \dim(\g')$. Define a map $F : \g'
  \to \g$ by $F(X_1') = X_1$, $F(Y_1') = Y_1$ and $\overline{F} =
  F|_{\qk'}$ is an isometry from $\qk'$ onto $\qk$. Then define a new
  Lie bracket on $\g$ by \[[X,Y]'' = F \left( [F^{-1}(X), F^{-1}(Y)]'
  \right), \ \forall X, Y \in \g.\] Denote by $(\g'', [\cdot,
  \cdot]'')$ this new Lie algebra. So $F$ is an i-isomorphism from
  $\g'$ onto $\g''$.

  Moreover $\g'' = (\CC X_1 \oplus \CC Y_1) \oplusp \qk$ is the double
  extension of $\qk$ by $\cb''$ with $\cb''= \overline{F} \ \cb'\
  \overline{F}^{-1}$. Then $\g$ and $\g'$ are i-isomorphic if, and
  only if, $\g$ and $\g''$ are i-isomorphism. Applying Proposition
  \ref{4.6}, this is the case if, and only if, there exists
  $\overline{A} \in \OO(\qk)$ such that $\cb'' = \lambda \
  \overline{A} \ \cb \ \overline{A}^{-1}$ for some non-zero complex
  $\lambda$. That implies \[ \cb' = \lambda \ ( \overline{F}^{-1} \
  \overline{A}) \ \cb \ ( \overline{F}^{-1} \overline{A})^{-1}\] and
  proves (1).

  \smallskip

\item We keep the notation in (1). We have that $\g$ and $\g'$ are
  isomorphic if, and only if, $\g$ and $\g''$ are isomorphic. Applying
  Proposition \ref{4.6}, $\g$ and $\g''$ are isomorphic if, and only
  if, there exists an invertible map $\overline{P} \in \Lc(\qk)$ and a
  non-zero $\lambda \in \CC$ such that $\cb'' = \lambda \ \overline{P}
  \ \cb \ \overline{P}^{-1}$ and $\overline{P}^* \ \overline{P} \ \cb
  = \cb$ and we conclude that $\cb' = \lambda \ \overline{Q} \ \cb \
  \overline{Q}^{-1}$ with $\overline{Q} = \overline{F}^{-1} \
  \overline{P}$. Finally, $\overline{F}^{-1} = \overline{Q} \
  \overline{P}^{-1}$ is an isometry from $\qk$ to $\qk'$.

  On the other hand, if $\cb' = \lambda \ \overline{Q} \ \cb \
  \overline{Q}^{-1}$ and $\overline{P}^* \ \overline{P} \ \cb = \cb$
  with $\overline{P} = \overline{F} \ \overline{Q}$ for some isometry
  $\overline{F} : \qk' \to \qk$, then construct $\g''$ as in (1). We
  deduce $\cb'' = \lambda \ \overline{P} \ \cb \ \overline{P}^{-1}$
  and $\overline{P}^* \ \overline{P} \ \cb = \cb$. So, by Proposition
  \ref{4.6}, $\g$ and $\g''$ are isomorphic and therefore, $\g$ and
  $\g'$ are isomorphic. 

\smallskip

\end{enumerate}

\end{proof}

\begin{rem} \label{4.8} Let $\g$ be a solvable singular quadratic Lie
  algebra. Consider $\g$ as a double extension of two quadratic
  vectors spaces $\qk$ and $\qk'$: \[ \g = (\CC X_1 \oplus \CC Y_1)
  \oplusp \qk \ \text{ and } \ \g = (\CC X_1' \oplus \CC Y_1') \oplusp
  \qk'.\] Let $\cb = \ad(Y_1)|_\qk$ and $\cb' =\ad(Y_1')|_{\qk'}$
  Since $\Id_\g$ is obviously an i-isomorphism, there exists an
  isometry $\overline{A} : \qk \to \qk'$ and a non-zero $\lambda \in
  \CC$ such that \[ \cb' = \lambda \ \overline{A} \ \cb \
  \overline{A}^{-1}.\]
\end{rem}

\begin{rem}
  A weak form of Corollary \ref{4.7} (1) was stated in \cite{FS}, in
  the case of i-isomorphisms satisfying some (dispensable)
  conditions. So (1) is an improvement. To our knowledge, (2) is
  completely new. Corollary \ref{4.7} and Remark \ref{4.8} can be
  applied directly to solvable singular Lie algebras: by Propositions
  \ref{4.4} and \ref{4.5}, they are double extensions of quadratic
  vector spaces by skew-symmetric maps.
\end{rem}

\subsection{} We shall now classify solvable singular Lie algebra
structures on $\CC^{n+2}$ up to i-isomorphism in terms of
$\OO(n)$-orbits in $\ps(\ok(n))$. We need the Lemma below:

\begin{lem} \label{4.9} Let $V$ be a quadratic vector space such that
  $V = (\CC X_1 \oplus \CC Y_1) \oplusp \qk'$ with $X_1$, $Y_1$
  isotropic elements and $B(X_1, Y_1) = 1$. Let $\g$ be a solvable
  singular quadratic Lie algebra with $\dim(\g) = \dim(V)$. Then,
  there exists a skew-symmetric map $\cpb : \qk' \to \qk'$ such that
  $V$ considered as the the double extension of $\qk'$ by $\cpb$ is
  i-isomorphic to $\g$.
\end{lem}

\begin{proof}
  By Propositions \ref{4.4} and \ref{4.5}, $\g$ is a double
  extension. Let us write $\g = (\CC X_0 \oplus \CC Y_0) \oplusp \qk$
  and $\cb = \ad(Y_0)|_\qk$. Define $A : \g \to V$ by $A(X_0) = X_1$,
  $A(Y_0) = Y_1$ and $\overline{A} = A |_\qk$ any isometry from $\qk
  \to \qk'$. It is clear that $A$ is an isometry from $\g$ to
  $V$. Now, define the Lie bracket on $V$ by: \[ [X,Y] = A \left( [
    A^{-1} (X), A^{-1}(Y)] \right), \ \forall \ X, Y \in V.\] Then $V$
  is a quadratic Lie algebra, that is i-isomorphic to $\g$, by
  definition. Moreover, $V$ is obviously a double extension of $\qk'$
  by $\cpb = \overline{A} \ \cb \ \overline{A}^{-1}$.
\end{proof}

We can now apply our results to the classification (up to
i-isomorphism) of solvable elements of $\Ss(n+2)$ (the set of singular
Lie algebras structures on $\CC^{n+2}$), for $n \geq 2$. We denote by
$\Sso(n+2)$ the set of solvable elements of $\Ss(n+2)$. Given $\g \in
\Ss(n+2)$, we denote by $[\g]_{\mathrm{i}}$ its i-isomorphism class
and by $\Ssoi(n+2)$ the set of classes. For $[\cb] \in \ps(\ok(n))$,
we denote by $\Os_{[\cb]}$ its $\OO(n)$-adjoint orbit and by
$\widetilde{\ps(\ok(n))}$ the set of orbits.

\begin{prop}\label{4.11}
  There exists a bijection $\theta : \widetilde{\ps(\ok(n))} \to
  \Ssoi(n+2)$.
\end{prop}

\begin{proof}
  We consider $\Os_{[\cb]} \in \widetilde{\ps(\ok(n))}$. There is a
  double extension $\g$ of \linebreak $\qk =\spa \{E_2, \dots,
  E_{n+1}\}$ by $\cb$ realized on $\CC^{n+2} = (\CC E_1 \oplus \CC
  E_{n+2} ) \oplusp \qk$. Then, by Corollary \ref{4.7}, $\g \in \Sso
  (n+2)$ and $[\g]_{\mathrm{i}}$ does not depend on the choice of
  $\cb$. We define $\theta(\Os_{[\cb]}) = [\g]_\irm$. If $\g' \in
  \Sso(n+2)$ then by Lemma \ref{4.9}, $\g'$ can be realized (up to
  i-isomorphism) as a double extension on $\CC^{n+2} = (\CC E_1 \oplus
  \CC E_{n+2} ) \oplusp \qk$. So $\theta$ is onto. Finally, $\theta$
  is one-to-one by Corollary \ref{4.7}.
\end{proof}

\section{Nilpotent and diagonalizable cases}

\subsection{} Let us denote by $\Ns(n+2)$ the set of nilpotent
elements of $\Ss(n+2)$, for $n \geq 1$.  Given $\g \in \Ns(n+2)$, we
denote by $[\g]$ its isomorphism class and by $[\g]_\irm$ its
i-isomorphism class. The set $\Nsh(n+2)$ is the set of all isomorphism
classes and $\Nsi(n+2)$ is the set of all i-isomorphism classes of
elements in $\Ns(n+2)$.

Let $\Nc(n)$ be the set of non-zero nilpotent elements of
$\ok(n)$. Given $\cb \in \Nc(n)$, we denote by $\Oc_\cb$ its
$\OO(n)$-adjoint orbit. The set of nilpotent orbits is denoted by
$\widetilde{\Nc}(n)$.

\begin{lem} \label{5.1}
  Let $\cb $ and $\cb' \in \Nc(n)$. Then $\cb$ is conjugate to
  $\lambda \cb'$ modulo $\OO(n)$ for some non-zero $\lambda \in \CC$
  if, and only if, $\cb$ is conjugate to $\cb'$.
\end{lem}

\begin{proof}
  It is enough to show that $\cb$ and $\lambda \cb$ are conjugate, for
  any non-zero $\lambda \in \CC$. By \cite{CM}, there exists a
  $\slk(2)$-triple $\{X, H, \cb \}$ in $\ok(n)$ such that $[H,\cb] = 2
  \cb$, so $\exx^{t \ad(H)} (\cb) = \exx^{2t} \cb$, $\forall t \in
  \CC$. We choose $t$ such that $\exx^{2t} = \lambda$, then $\exx^{tH}
  \cb \exx^{-t H} = \lambda \cb$ and $\exx^{tH} \in \OO(n)$.
\end{proof}

\begin{prop} \label{5.2} One has:

\begin{enumerate}

\item Let $\g$ and $\g' \in \Ns(n+2)$. Then $\g$ and $\g'$ are
  isomorphic if, and only if, they are i-isomorphic, so $[\g]_\irm =
  [\g]$ and $\Nsi(n+2) = \Nsh(n+2)$.

\smallskip

\item There is a bijection $\tau : \widetilde{\Nc}(n) \to \Nsh(n+2)$.

\smallskip

\item $\Nsh(n+2)$ is finite.

\end{enumerate}

\end{prop}

\begin{proof} \hfill

\begin{enumerate}

\item Using Lemma \ref{4.9}, Proposition \ref{4.4}(3) and Corollary
  \ref{4.7}, it is enough to show that for $\cb$ and $\cb' \in
  \Nc(n+2)$, if there exists $P \in \GL(n)$ such that $\cb' = \lambda
  P \cb P^{-1}$, for some non-zero $\lambda \in \CC$, then $\cb$ and
  $\cb'$ are conjugate under $\OO(n)$. By Lemma \ref{5.1}, we can
  assume that $\lambda = 1$, and then the result is well known (see
  e.g. \cite{CM}).

\smallskip

\item As in the proof of Proposition \ref{4.11}, for a given $\Oc_\cb
  \in \widetilde{\Nc}(n)$, we construct the double extension $\g$ of
  $\qk =\spa \{E_2, \dots, E_{n+1}\}$ by $\cb$ realized on
  $\CC^{n+2}$. Then, by Proposition \ref{4.4} (3), $\g \in \Ns(n+2)$
  and $[\g]$ does not depend on the choice of $\cb$. We define
  $\tau(\Os_{\cb}) = [\g]$. Then by (1) and Corollary \ref{4.7},
  $\tau$ is one-to-one and onto.

\item $\Nsh(n+2)$ is finite since the set of nilpotent orbits
  $\widetilde{\Nc}(n)$ is finite (see e.g. \cite{CM}).

\end{enumerate}

\end{proof}

\begin{defn} \label{5.2bis}

  Let $p \in \NN \setminus \{ 0 \}$. We denote the {\em Jordan block
    of size} $p$ by $J_1 := (0)$ and for $p \geq 2$, \[J_p :
  = \begin{pmatrix} 0 & 1 & 0 & \dots & 0 \\ 0 & 0 & 1 & \dots & 0\\
    \vdots & \vdots & \dots & \ddots & \vdots \\ 0 & 0 & \dots & 0 & 1
    \\ 0 & 0 & 0 & \dots & 0 \end{pmatrix}.\]

Next, we define nilpotent Jordan-type Lie algebras. There are two
types:

\begin{itemize}

\item for $p \geq 2$, we consider $\qk = \CC^{2p}$ equipped with its
  canonical bilinear form $\overline{B}$ and the map $\cb_{2p}^J$ with
  matrix \[ \begin{pmatrix} J_p & 0 \\ 0 & -{}^t J_p\end{pmatrix}\] in
  the canonical basis. Then $\cb_{2p}^J \in \ok(2p)$ and we denote by
  $\jk_{2p}$ the double extension of $\qk $ by $\cb_{2p}^J$. So
  $\jk_{2p} \in \Ns(2p+2)$.

\item for $p \geq 1$, we consider $\qk = \CC^{2p+1}$ equipped with its
  canonical bilinear form $\overline{B}$ and the map $\cb_{2p+1}^J$
  with matrix \[ \begin{pmatrix} J_{p+1} & M \\ 0 & -{}^t
    J_p\end{pmatrix}\] in the canonical basis, where $M=(m_{ij})$
  denotes the $(p+1) \times p$-matrix with $m_{p+1,p} = -1$ and
  $m_{ij} = 0$ otherwise. Then $\cb_{2p+1}^J \in \ok(2p+1)$ and we
  denote by $\jk_{2p+1}$ the double extension of $\qk $ by
  $\cb_{2p+1}^J$. So $\jk_{2p+1} \in \Ns(2p+3)$.

\end{itemize}

\medskip

Lie algebras $\jk_{2p}$ or $\jk_{2p+1}$ will be called {\em nilpotent
  Jordan-type Lie algebras}.

\end{defn}

Let $n \in \NN$, $n \neq 0$. We consider partitions $[d] := (d_1,
\dots, d_r)$ of $n$ of a special type:

\begin{itemize}

\item each even $d_i$ must occur with even multiplicity.

\item $[d]$ can be written as $(p_1, p_1, p_2, p_2, \dots, p_k, p_k, 2
  q_1 +1, \dots 2 q_\ell +1)$ with all $p_i$ even, $p_1 \geq p_2 \geq
  \dots \geq p_k$ and $q_1 \geq q_2 \geq \dots \geq q_\ell$.

\end{itemize}

We denote by $\Pc'(n)$ the set of partitions satisfying the above
conditions. To every $[d] \in \Pc'(n)$, we associate a map $\cb_{[d]}
\in \ok(n)$: write \linebreak $[d] = (p_1, p_1, p_2, p_2, \dots, p_k,
p_k, 2 q_1 +1, \dots 2 q_\ell +1)$. Then $\cb_{[d]}$ is the map with
matrix
\[ \diag_{2k + \ell}(\cb^J_{2p_1}, \cb^J_{2p_2}, \dots, \cb^J_{2p_k},
\cb^J_{2q_1+1}, \dots, \cb^J_{2q_\ell +1} ).\] in the canonical basis
of $\CC^n$.

Moreover, we denote by $\g_{[d]}$ the double extension of $\CC^n$ by
$\cb_{[d]}$. Then $\g_{[d]} \in \Ns(n+2)$ and $\g_{[d]}$ is an
amalgamated product of nilpotent Jordan-type Lie algebras, more
precisely, \[ \g_{[d]} = \jk_{2p_1} \ta \jk_{2p_2} \ta \dots \ta
\jk_{2p_k} \ta \jk_{2q_1 + 1} \ta \dots \ta \jk_{2q_\ell +1}.\] The
following fundamental result classifies all nilpotent $\OO(n)$-orbits
in $\ok(n)$ (see \cite{CM}).

\begin{lem} \label{5.3}
  The map $[d] \mapsto \cb_{[d]}$ from $\Pc'(n)$ to $\ok(n)$ induces a
  bijection from $\Pc'(n)$ onto $\widetilde{\Nc}(n)$.
\end{lem}

Using Propositions \ref{5.2} and \ref{5.3}, we deduce:

\begin{prop} \label{5.4} \hfill

\begin{enumerate}

\item The map $[d] \mapsto \g_{[d]}$ from $\Pc'(n)$ to $\Ns(n+2)$
  induces a bijection from $\Pc'(n)$ onto $\Nsh(n+2)$.

\smallskip

\item Each nilpotent singular $n+2$-dimensional Lie algebra is
  i-isomorphic to a unique amalgamated product $\g_{[d]}$, $[d] \in
  \Pc'(n)$ of nilpotent Jordan-type Lie algebras.

\smallskip

\end{enumerate}
\end{prop}

\subsection{} We introduce some notation:

\begin{defn}
  Let $\g$ be a solvable singular quadratic Lie algebra and write $\g
  = (\CC X_0 \oplus \CC Y_0) \oplusp \qk$ a decomposition of $\g$ as a
  double extension (Proposition \ref{4.4} and Lemma \ref{4.5}). Let
  $\cb = \ad(Y_0)|_\qk$. We say that $\g$ is a {\em diagonalizable} if
  $\cb$ is diagonalizable.

  We denote by $\Ds(n+2)$ the set of such structures on the quadratic
  space $\CC^{n+2}$, by $\Ds_{\mathrm{red}}(n+2)$ the reduced ones, by
  $\Dsh(n+2)$, $\Dsi(n+2)$, $\Dsh_{\mathrm{red}}(n+2)$,
  $\Dsi_{\mathrm{red}}(n+2)$ the corresponding sets of isomorphism and
  i-isomorphism classes of elements in $\Ds(n+2)$ and
  $\Ds_{\mathrm{red}}(n+2)$.
\end{defn}

Remark that the property of being diagonalizable does not depend on
the chosen decomposition of $\g$ (see Remark \ref{4.8}). By Corollary
\ref{4.7} and a proof completely similar to Proposition \ref{4.11} or
Proposition \ref{5.2}, we conclude:

\begin{prop} \label{5.5} There is a bijection between $\Dsi(n+2)$ and
  the set of semi-simple $\OO(n)$-orbits in $\ps(\ok(n))$. The same
  result holds for $\Dsi_{\mathrm{red}}(n+2)$ and semi-simple
  invertible orbits in $\ps(\ok(n))$.
\end{prop}

\begin{proof}
  Proceed exactly as in Proposition \ref{4.11} or Proposition
  \ref{5.2}, but notice that a a diagonalizable $\cb$ satisfies
  $\ker(\cb) \subset \im(\cb)$ if, and only if, $\ker(\cb) = \{ 0\}$.
\end{proof}

\subsection{} \label{5.4s} The classification of semi-simple adjoint
orbits of a semi-simple Lie algebra $\g$ is fully known (see
e.g. \cite{CM}). Given a Cartan subalgebra $\hk$ of $\g$, there is a
bijection between the set of semi-simple adjoint orbits and $\hk / W$,
where $W$ is the Weyl group.

Here, we deal with $\OO(n)$-adjoint and not $\SO(n)$-adjoint
orbits. Hence, slight changes must be done. Let us recall the result:
write $n = 2p$ if $n$ is even and $n = 2p+1$ if $n$ is odd. Let $\hk$
be a Cartan subalgebra given by the vector space of diagonal matrices
of type $\diag_{2p}(\lambda_1, \dots, \lambda_p, \lambda_1, \dots,
\lambda_p)$ if $n$ is even and of type \linebreak
$\diag_{2p+1}(\lambda_1, \dots, \lambda_p, 0, -\lambda_1, \dots,
-\lambda_p)$ if $n$ is odd. Any diagonalizable $\cb \in \ok(n)$ is
conjugate to an element of $\hk$ (see Appendix 1 for a direct
proof). If $\cb$ is invertible, then $n$ is even (see Appendix 1).

If $n$ is even, the Weyl group consists of all permutations and even
sign changes of $(\lambda_1, \dots, \lambda_p)$. Thus, to describe
$\OO(n)$-orbits we must admit any number of sign changes. We denote by
$G_p$ the corresponding group. If $n$ is odd, the Weyl group is $G_p$
and there is nothing to add.

However, we are interested in $\OO(n)$-orbits in $\ps(\ok(n))$. So, we
must add maps $(\lambda_1, \dots, \lambda_p) \mapsto \lambda
(\lambda_1, \dots, \lambda_p)$, $\forall \lambda \in \CC$, $\lambda
\neq 0$ to the group $G_p$. We obtain a group denoted by $H_p$. Now,
let $\Lambda_p = \{ (\lambda_1, \dots, \lambda_p) \mid \lambda_1,
\dots, \lambda_p \in \CC, \lambda_i \neq 0 \ \text{ for some } \ i \}$
and $\Lambda_p^+ = \{ (\lambda_1, \dots, \lambda_p) \mid \lambda_1,
\dots, \lambda_p \in \CC, \lambda_i \neq 0, \forall \ i \}$.

By Proposition \ref{5.5}, we obtain the Corollary:

\begin{cor}\label{5.8}
  There is a bijection between $\Dsi(n+2)$ and $\Lambda_p /
  H_p$. Moreover, if $n = 2p +1$, $\Dsh_{\mathrm{red}}(n+2) =
  \emptyset$ and if $n =2p$, then $\Dsh_{\mathrm{red}}(2p+2)$ is in
  bijection with $\Lambda_p^+/H_p$.
\end{cor}

\subsection{} To go further in the study of diagonalizable reduced
case, we need the following Lemma that will also be used in Section 6:

\begin{lem} \label{5.7} \hfill 

  Let $\g'$ and $\g''$ be solvable singular quadratic Lie algebras,
  $\g' = (\CC X_1' \oplus \CC Y_1') \oplusp \qk'$ a decomposition of
  $\g'$ as a double extension and $\cb' = \ad(Y_1')|_{\qk'}$. We
  assume that $\cb'$ is invertible. Then $\g'$ and $\g''$ are
  isomorphic if, and only if, they are i-isomorphic.
\end{lem}

\begin{proof}
Write $\g'' = (\CC X_1'' \oplus \CC Y_1'') \oplusp \qk''$ a decomposition of $\g''$
  as a double extension and $\cb'' = \ad(Y_1'')|_{\qk''}$.

  Assume that $\g'$ and $\g''$ are isomorphic. By Corollary \ref{4.7},
  there exist $\overline{Q} : \qk' \to \qk''$ and $\overline{P} \in
  \Lc(\qk')$ such that $\overline{Q} \ \overline{P}^{-1}$ is an
  isometry, $\overline{P}^* \ \overline{P} \ \cb' = \cb'$ and $\cb'' =
  \lambda \ \overline{Q} \ \cb' \ \overline{Q}^{-1}$ for some non-zero
  $\lambda \in \CC$. But $\cb'$ is invertible, so $\overline{P}^* \
  \overline{P} = \Id_{\qk'}$. Therefore, $\overline{P}$ is an isometry
  of $\qk'$ and then $\overline{Q}$ is an isometry from $\qk'$ to
  $\qk''$. The conditions of Corollary \ref{4.7} (1) are satisfied, so
  $\g'$ and $\g''$ are i-isomorphic.
\end{proof}

\begin{cor} \label{5.9} One has: \[ \Dsh_{\mathrm{red}} (2p+2) =
  \Dsi_{\mathrm{red}} (2p+2), \ \forall \ p \geq 1.\]
\end{cor}

Next, we describe diagonalizable reduced singular Lie algebras using
amalgamated products. First, let $\g_4(\lambda)$ be the double
extension of $\qk = \CC^2$ by $\cb = \begin{pmatrix} \lambda & 0 \\ 0
  & - \lambda \end{pmatrix}$. By Lemma \ref{4.3}, $\g_4(\lambda)$ is
i-isomorphic to $\g_4(1)$, call it $\g_4$.

\begin{prop}\label{5.10}
  Let $(\g,B)$ be a diagonalizable reduced singular Lie algebra. Then
  $\g$ is an amalgamated product of singular Lie algebras all
  i-isomorphic to $\g_4$.
\end{prop}

\begin{proof}
  We write $\g = (\CC X_0 \oplus \CC Y_0) \oplusp \qk$, $C =
  \ad(Y_0)$, $\cb = C|_\qk$ and $\overline{B} = B_{\qk \times
    \qk}$. Then $\cb$ is a diagonalizable invertible element of
  $\ok(\qk, \overline{B})$. Apply Appendix 1 to obtain a basis $\{e_1,
  \dots, e_p, f_1, \dots, f_p \}$ of $\qk$ and $\lambda_1, \dots,
  \lambda_p \in \CC$, all non-zero, such that $B(e_i, e_j) = B(f_i;
  f_j) = 0$, $B(e_i,f_j) = \delta_{ij}$ and $\cb(e_i) = \lambda_i
  e_i$, $\cb(f_i) = - \lambda_i f_i$, for all $1 \leq i,j \leq p$. Let
  $\qk_i = \spa \{e_i, f_i \}$, $1 \leq i \leq p$. Then \[\qk =
  \oplusp_{i = 1}^p \qk_i.\] Furthermore, $\hk_i = (\CC X_0 \oplus \CC
  Y_0) \oplusp \qk_i$ is a Lie subalgebra of $\g$ for all $1 \leq i
  \leq p$ and \[ \g = \hk_1 \ta \hk_2 \ta \dots \ta \hk_p \ \text{
    with } \ \hk_i \iiso \g_4(\lambda_i) \iiso \g_4.\]
\end{proof}

\begin{rem}\label{5.12}
  For non-zero $\lambda$, $\mu \in \CC$, consider the amalgamated
  product: \[ \g(\lambda, \mu) = \g_4(\lambda) \ta \g_4(\mu).\] Then
  $\g(\lambda, \mu)$ is the double extension of $\CC^4$
  by \[ \begin{pmatrix} \lambda & 0 & 0 & 0 \\ 0 & \mu & 0 & 0 \\ 0 &
    0 & - \lambda & 0 \\ 0 & 0 & 0 & - \mu \end{pmatrix}.\] Therefore
  $\g(\lambda, \mu)$ is isomorphic to $\g(1,1)$ if, and only if, $\mu
  = \pm \lambda$ (Lemma \ref{5.7} and (\ref{5.4s})). So, though
  $\g_4(\lambda)$ and $\g_4(\mu)$ are i-isomorphic to $\g_4$, the
  amalgamated product $\g(\lambda, \mu)$ is not even isomorphic to
  $\g(1,1) = \g_4 \ta \g_4$ if $\mu \neq \pm \lambda$. This
  illustrates that amalgamated products may have a rather bad behavior
  with respect to isomorphisms.
\end{rem}

\section{The general case}

\subsection{} Let $\g$ be a solvable singular quadratic Lie
algebra. We fix a realization of $\g$ as a double extension, $\g =
(\CC X_0 \oplus \CC Y_0) \oplusp \qk$ (Propositions \ref{4.4} and
{4.5}). Let $C = \ad(Y_0)$ and $\cb = C|_\qk$. We consider the Fitting
decomposition of $\cb$: \[ \qk = \qk_N \oplus \qk_I, \] where $\qk_N$
ans $\qk_I$ are $\cb$-stable, $\cb_N = \cb|_{\qk_N}$ is nilpotent and
$\cb_I = \cb|_{\qk_I}$ is invertible.

Since $\cb$ is skew-symmetric, one has $\qk_I =
\qk_N^\perp$. Therefore, the restrictions $\overline{B}_N =
\overline{B}|_{\qk_N \times \qk_N}$ and $\overline{B}_I =
\overline{B}|_{\qk_I \times \qk_I}$ of $\overline{B} = B|_{\qk \times
  \qk}$ are non degenerate, $\cb_N$ and $\cb_I$ are skew-symmetric and
$[\qk_I, \qk_N] = 0$. Let $\g_N = (\CC X_0 \oplus \CC Y_0) \oplusp
\qk_N$ and $\g_I = (\CC X_0 \oplus \CC Y_0) \oplusp \qk_I$. Then
$\g_N$ and $\g_I$ are Lie subalgebras of $\g$, $\g_N$ is the double
extension of $\qk_N$ by $\cb_N$, $\g_I$ is the double extension of
$\qk_I$ by $\cb_I$ and $\g_N$ is a nilpotent singular quadratic Lie
algebra. To study $\g_I$ , we introduce the following definition:

\begin{defn}
  A double extension is called an {\em invertible quadratic Lie
    algebra} if the corresponding skew-symmetric map is invertible.
\end{defn}

\begin{rem}\label{6.2} \hfill

\begin{itemize}

\item By Remark \ref{4.8}, the property of being an invertible quadratic Lie
algebra does not depend on the chosen decomposition. 

\item By Appendix 1, the dimension of an invertible quadratic Lie
  algebra is even.

\item By Lemma \ref{5.7}, two invertible quadratic Lie algebras are
  isomorphic if, and only if, they are i-isomorphic.
\end{itemize}

\end{rem}

With the above definition, $\g_I$ is an invertible quadratic Lie
algebra and we have \[\g = \g_N \ta \g_I.\]

\begin{defn}
  The Lie subalgebras $\g_N$ and $\g_I$ are respectively the {\em
    nilpotent} and {\em invertible Fitting components} of $\g$.
\end{defn}

This definition is justified by:

\begin{prop}\label{6.4}
  Let $\g$ and $\g'$ be solvable singular quadratic Lie algebras and
  $\g_N$, $\g_I$, $\g'_N$, $\g'_I$ be their Fitting components. Then

\begin{enumerate}

\item $\g \iiso \g'$ if, and only if, $\g_N \iiso \g'_N$ and $\g_I
  \iiso \g_I'$. The result remains valid if we replace $\iiso$ by
  $\simeq$.

\smallskip

\item $\g \simeq \g'$ if, and only of, $\g \iiso \g'$.

\end{enumerate}

\end{prop}

\begin{proof}
  We assume that $\g \simeq \g'$. Then by Corollary \ref{4.7}, there
  exists an invertible $\overline{P} : \qk \to \qk'$ and a non-zero
  $\lambda \in \CC$ such that $\cb' = \lambda \ \overline{P} \ \cb \
  \overline{P}^{-1}$, so $\qk_N' = \overline{P} (\qk_N)$ and $\qk_I' =
  \overline{P} (\qk_I)$, then $\dim(\qk'_N) = \dim(\qk_N)$ and
  $\dim(\qk'_I) = \dim(\qk_I)$. Thus, there exist isometries $F_N :
  \qk'_N \to \qk_N$ and $F_I : \qk'_I \to \qk_I$ and we can define an
  isometry $\overline{F} : \qk' \to \qk$ by $\overline{F}(X_N' + X_I')
  = F_N(X_N') + F_I(X_I')$, $\forall X_N' \in \qk_N'$ and $ X_I' \in
  \qk_I'$. We now define $F : \g' \to \g$ by $F(X_1') = X_1$, $F(Y_1')
  = Y_1$, $F|_{\qk'} = \overline{F}$ and a new Lie bracket on $\g$
  : \[ [X,Y]'' = F\left( [F^{-1}(X), F^{-1}(Y)]' \right), \ \forall X,
  Y \in \g.\]

  Call $\g''$ this new quadratic Lie algebra. We have $\g'' = ( \CC
  X_1 \oplus \CC Y_1) \oplusp \qk$, i.e., $\qk'' = \qk$ and $\cb'' =
  \overline{F} \ \cb' \ \overline{F}^{-1}$. So $\qk_N'' = F(\qk'_N) =
  \qk_N$ and $\qk_I'' = F(\qk'_I) = \qk_I$. But $\g \simeq \g''$, so
  there exists an invertible $Q : \qk \to \qk$ such that $\cb'' =
  \lambda \ \overline{Q} \ \cb \ \overline{Q}^{-1}$ for some non-zero
  $\lambda \in \CC$ (Corollary \ref{4.7}). It follows that $\qk_N'' =
  Q(\qk_N)$ and $\qk_I = Q(\qk_I)$, so $Q (\qk_N) = \qk_N$ and $Q
  (\qk_I) = \qk_I$.

  Moreover, we have $Q^* \ Q \ \cb = \cb$ (Corollary \ref{4.7}), so
  $Q^* \ Q \ \cb^k = \cb^k$ for all $k$. There exists $k$ such that
  $\qk_I = \im(C^k)$ and $(Q^* \ Q \ \cb^k)(X) = \cb^k(X)$, for all $X
  \in \g$. So $Q^* Q|_{\qk_I} = \Id_{\qk_I}$ and $Q_I = Q|_{\qk_I}$ is
  an isometry. Since $\cb_I'' = \lambda \ \overline{Q}_I \ \cb_I \
  \overline{Q}_I^{-1}$, then $\g_I \iiso \g_I''$ (Corollary
  \ref{4.7}).
 
  Let $Q_N = Q|_{\qk_N}$. Then $\cb_N'' = \lambda \ \overline{Q}_N \
  \cb_N \ \overline{Q}_N^{-1}$ and $Q_N^* \ Q_N \ \cb_N = \cb_N$, so
  by Corollary \ref{4.7}, $\g_N \simeq \g_N''$. Since $\g_N$ and
  $\g_N''$ are nilpotent, then $\g_N '' \iiso \g_N$ by Proposition
  \ref{5.2}.

  Conversely, assume that $\g_N \simeq \g'_N$ and $\g_I \simeq
  \g_I'$. Then $\g_N \iiso \g'_N$ and $\g_I \iiso \g_I'$ by
  Proposition \ref{5.2} and Lemma \ref{5.7}.

  So, there exist isometries $P_N : \g_N \to \g_N'$, $P_I : \g_I \to
  \g_I'$ and non-zero $\lambda_N$ and $\lambda_I \in \CC$ such that
  $\cb_N' = \lambda_N \ \overline{P}_N \ \cb_N \ \overline{P}_N^{-1}$
  and $\cb_I' = \lambda_I \ \overline{P}_I \ \cb_I \
  \overline{P}_I^{-1}$. By Lemma \ref{5.1}, since $\g_N$ and $\g_N'$
  are nilpotent, we can assume that $\lambda_N = \lambda_I =
  \lambda$. Now we define $P : \qk \to \qk'$ by $P(X_N + X_I) =
  P_N(X_N) + P_I(X_I)$, $\forall X_N \in \qk_N$, $X_I \in \qk_I$, so
  $P$ is an isometry. Moreover, since $\cb(X_N + X_I) = \cb_N(X_N) +
  \cb_I(X_I)$, $\forall X_N \in \qk_N$, $X_I \in \qk_I$ and $\cb'(X_N'
  + X_I') = \cb_N'(X_N') + \cb_I'(X_I')$, for all $X_N' \in \qk_N$,
  $X_I' \in \qk_I$, we conclude $\cb' = \lambda \ P \ \cb P^{-1}$ and
  finally, $\g \iiso \g'$, by Corollary \ref{4.7}.
\end{proof}

\begin{rem}
  The class of solvable singular quadratic Lie algebras has the
  remarkable property that two Lie algebras in this class are
  isomorphic if, and only if, they are i-isomorphic. In addition, the
  Fitting components do not depend on the realizations of the Lie
  algebra as a double extension and they completely characterize the
  Lie algebra (up to isomorphism).
\end{rem}

\subsection{} To classify all solvable singular Lie algebras (up to
isomorphism), we have to classify nilpotent and invertible ones (see
Proposition \ref{6.4}). The nilpotent case is completely achieved in
Proposition \ref{5.4}, so we are left with the invertible case.

For $p \geq 1$ and $\lambda \in \CC$, let $J_p(\lambda) =
\diag_p(\lambda, \dots, \lambda) + J_p$ and \[\cb_{2p}^J (\lambda)
= \begin{pmatrix} J_p(\lambda) & 0 \\ 0 & - {}^t
  J_p(\lambda) \end{pmatrix}.\] Then $\cb_{2p}^J(\lambda) \in
\ok(2p)$.

\begin{defn}
  For $\lambda \in \CC$, let $\jk_{2p}(\lambda)$ be the double
  extension of $\CC^{2p}$ by $\cb_{2p}^J(\lambda)$. We say that
  $\jk_{2p}(\lambda)$ is a {\em Jordan-type quadratic Lie algebra}.

  When $\lambda = 0$ and $p \geq 2$, we recover the nilpotent
  Jordan-type Lie algebras $\jk_{2p}$ from Definition \ref{5.2bis}.

  When $\lambda \neq 0$, $\jk_{2p}(\lambda)$ is an invertible singular
  quadratic lie algebra and \[\jk_{2p}(- \lambda) \simeq
  \jk_{2p}(\lambda).\]
\end{defn}

\begin{prop}\label{6.5}
  Let $\g$ be a solvable singular quadratic Lie algebra. Then $\g$ is
  an invertible Lie algebra if, and only if, $\g$ is an amalgamated
  product of Lie algebras all i-isomorphic to Jordan-type Lie algebras
  $\jk_{2p} (\lambda)$, with $\lambda \neq 0$.
\end{prop}

\begin{proof}
  Let $\g = (\CC X_0 \oplus \CC Y_0) \oplusp \qk$, $B$ be the bilinear
  form of $\g$, $\overline{B} = B|_{\qk \times \qk}$, $C = \ad(Y_0)$
  and $\cb = C|_\qk \in \ok(\qk, \overline{B})$. We decompose $\cb$
  into its semi-simple and nilpotent parts, $\cb = \overline{S} +
  \overline{N}$. It is well known that $\overline{S}$ and
  $\overline{N} \in \ok(\qk, \overline{B})$.

  Let $\Lambda \subset \CC \setminus \{0\}$ be the spectrum of
  $\overline{S}$. We have that $\lambda \in \Lambda$ if, and only if,
  $-\lambda \in \Lambda$ (see Appendix 1). Let $V_\lambda$ be the
  eigenspace corresponding to the eigenvalue $\lambda$. We have
  $\dim(V_\lambda) = \dim(V_{-\lambda})$. Denote by $\qk(\lambda)$ the
  direct sum $\qk(\lambda) = V_\lambda \oplus V_{-\lambda}$. If $\mu
  \in \Lambda$, $\mu \neq \pm \lambda$, then $\qk(\lambda)$ and
  $\qk(\mu)$ are orthogonal (Appendix 1). Choose $\Lambda_+$ such that
  $\Lambda = \Lambda_+ \cup \left( - \Lambda_+ \right)$ and $
  \Lambda_+ \cap \left( - \Lambda_+ \right) = \emptyset$. We have (see
  Appendix 1): \[ \qk = \oplusp_{\lambda \in \Lambda_+}
  \qk(\lambda).\] So the restriction $B_\lambda = B|_{\qk(\lambda)
    \times \qk(\lambda)}$ is non degenerate. Moreover, $V_\lambda$ and
  $V_{-\lambda}$ are maximal isotropic subspaces in $\qk(\lambda)$.

  Now, consider the map $\Psi : V_{-\lambda} \to V_\lambda^*$ defined
  by $\Psi(u)(v) = B_\lambda(u,v)$, $\forall u \in V_{-\lambda}$, $v
  \in V_{\lambda}$. Then $\Psi$ is an isomorphism. Given any basis
  $\Bc(\lambda) = \{ e_1(\lambda), \dots, e_{n_\lambda}(\lambda) \}$
  of $V_\lambda$, there is a basis $\Bc(-\lambda) = \{ e_1(-\lambda),
  \dots, e_{n_\lambda}(-\lambda) \}$ of $V_{-\lambda}$ such that
  \linebreak $B_\lambda(e_i(\lambda), e_j(-\lambda)) = \delta_{ij}$,
  $\forall 1 \leq i, j \leq n_\lambda$: simply define $e_i(-\lambda) =
  \psi^{-1}(e_i(\lambda)^*)$, for all $1 \leq i \leq n_\lambda$.

  Remark that $\overline{N}$ and $\overline{S}$ commute, so
  $\overline{N}(V_\lambda) \subset V_\lambda$, $\forall \lambda \in
  \Lambda$. Define $\overline{N}_\lambda =
  \overline{N}|_{\qk(\lambda)}$, then $\overline{N}_\lambda \in
  \ok(\qk(\lambda), B_\lambda)$. Hence, if
  $\overline{N}_\lambda|_{V_\lambda}$ has a matrix $M_\lambda$ with
  respect to $\Bc(\lambda)$, then
  $\overline{N}_\lambda|_{V_{-\lambda}}$ has a matrix $-{}^tM_\lambda$
  with respect to $\Bc(-\lambda)$. We choose the basis $\Bc(\lambda)$
  such that $M_\lambda$ is of Jordan type, i.e. \[\Bc(\lambda) =
  \Bc(\lambda, 1) \cup \dots \cup \Bc(\lambda, r_\lambda),\] the
  multiplicity $m_\lambda$ of $\lambda$ is $m_\lambda =
  \sum_{i=1}^{r_\lambda} d_\lambda(i)$ where $d_\lambda(i) = \sharp
  \Bc(\lambda, i)$ and \[M_\lambda = \diag_{n_\lambda} \left(
    J_{d_\lambda(1)}, \dots, J_{d_\lambda(r_\lambda)} \right).\] The
  matrix of $C|_{\qk(\lambda)}$ written on the basis $\Bc(\lambda)
  \cup \Bc(-\lambda)$ is: \[ \diag_{n_\lambda} \left(
    J_{d_\lambda(1)}(\lambda), \dots,
    J_{d_\lambda(r_\lambda)}(\lambda), -{}^t
    J_{d_\lambda(1)}(\lambda), \dots, -{}^t
    J_{d_\lambda(r_\lambda)}(\lambda) \right).\]

  Let $\qk(\lambda,i)$ be the subspace generated by $\Bc(\lambda, i)
  \cup \Bc(-\lambda, i)$, for all $1 \leq i \leq r_\lambda$ and let
  $C(\lambda,i) = C|_{\qk(\lambda,i)}$. We have \[ \qk(\lambda) =
  \oplusp_{1 \leq i \leq r_\lambda} \qk(\lambda,i).\] The matrix of
  $C(\lambda,i)$ written on the basis of $\qk(\lambda, i)$ is $C_{2
    d_\lambda(i)}^J (\lambda)$. Let $\g(\lambda,i)$, $\lambda \in
  \Lambda_+$, $1 \leq i \leq r_\lambda$ be the double extension of
  $\qk(\lambda,i)$ by $C(\lambda,i)$. Then $\g(\lambda, i)$ is
  i-isomorphic to $\jk_{2 d_\lambda(i)}(\lambda)$. But \[ \qk =
  \underset{\underset{1 \leq i \leq r_\lambda}{\lambda \in
      \Lambda_+}}{\oplusp} \qk(\lambda,i) \ \text{ and } \
  C|_{\qk(\lambda,i)} = C(\lambda, i).\] Therefore, $\g$ is the
  amalgamated product \[ \g= \underset{\underset{1 \leq i \leq
      r_\lambda}{\lambda \in \Lambda_+}}{\ta} \g(\lambda,i).\]

\end{proof}

\subsection{} Denote by $\Si(2p+2)$ the set of invertible singular Lie
algebra structures on $\CC^{2p+2}$, by $\Sih(2p+2)$ the set of
isomorphism (or i-isomorphism) classes of $\Si(2p+2)$. Next, we will
give a classification of $\Si(2p+2)$. Using Propositions \ref{6.4} and
\ref{5.4}, a classification of $\Ssoh(n+2)$ can finally be achieved.

We shall need the following Lemma;

\begin{lem} \label{6.6} Let $(V,B)$ be a quadratic vector space. We
  assume that $V = V_+ \oplus V_-$ with $V_\pm$ totally isotropic
  vector subspaces.

\begin{enumerate}

\item Let $N \in \Lc(V)$ such that $N(V_\pm) \subset V_\pm$. We define
  maps $N_\pm$ by $N_+|_{V_+} = N|_{V_+}$, $N_+|_{V_-} = 0$,
  $N_-|_{V_-} = N|_{V_-}$ and $N_-|_{V_+} = 0$. Then $N \in \ok(V)$
  if, and only if, $N_- = -N_+^*$ and, in this case, $N = N_+ -
  N_+^*$.

\smallskip

\item Let $U_+ \in \Lc(V)$ such that $U_+$ is invertible, $U_+(V_+) =
  V_+$ and $U_+|_{V_-} = \Id_{V_-}$. We define $U \in \Lc(V)$ by
  $U|_{V_+} = U_+$ and $U|_{V_-} = \left(U_+^{-1}\right)^*$. Then $U
  \in \OO(V)$.

\smallskip

\item Let $N' \in \ok(V)$ such that $N'$ satisfies the assumptions of
  (1). Define $N_\pm$ as in (1). Moreover, we assume that there exists
  $U_+ \in \Lc(V_+)$, $U_+$ invertible such that \[ N_+'|_{V_+} =
  \left( U_+ \ N_+ \ U_+^{-1} \right)|_{V_+}.\] We extend $U_+$ to $V$
  by $U_+|_{V_-} = \Id_{V_-}$ and define $U \in \OO(V)$ as in
  (2). Then \[ N' = U \ N \ U^{-1}.\]

\end{enumerate}

\end{lem}

\begin{proof}
  The proof is a straightforward computation.
\end{proof}

Let us now consider $C \in \ok(n)$, $C$ invertible. Then, $n$ is even,
$n = 2p$ (see Appendix 1). We decompose $C = S + N$ into semi-simple
and nilpotent parts, $S$, $N \in \ok(2p)$. We have $\lambda \in
\Lambda$ if, and only if, $-\lambda \in \Lambda$ (Appendix 1), where
$\Lambda$ is the spectrum of $C$. Also $m(\lambda) = m(-\lambda)$, for
all $\lambda \in \Lambda$ with multiplicity $m(\lambda)$. Since $N$
and $S$ commute, we have $N(V(\pm \lambda)) \subset V(\pm \lambda)$
where $V_\lambda$ is the eigenspace of $S$ corresponding to $\lambda
\in \Lambda$. Denote by $W(\lambda)$ the direct sum \[W(\lambda) =
V_\lambda \oplus V_{-\lambda}.\]

Define the equivalence relation $\Rs$ on $\Lambda$ by: \[ \lambda \Rs
\mu \ \text{ if, and only if, } \ \lambda = \pm \mu.\] Then \[\CC^{2p}
= \oplusp_{\lambda \in \Lambda / \Rs} W(\lambda),\] and each
$(W(\lambda), B_\lambda)$ is a quadratic vector space with $B_\lambda
= B|_{W(\lambda) \times W(\lambda)}$.

Fix $\lambda \in \Lambda$. We write $W(\lambda) = V_+ \oplus V_-$ with
$V_\pm = V_{\pm \lambda}$. Then, with the notation in Lemma \ref{6.6},
define $N_{\pm \lambda} = N_\pm$. Since $N|_{V_-} = - N_\lambda^*$, it
is easy to verify that the matrices of $N|_{V_+}$ and $N|_{V_-}$ have
the same Jordan form. Let $(d_1(\lambda), \dots,
d_{r_\lambda}(\lambda))$ be the size of the Jordan blocks in the
Jordan decomposition of $N|_{V_+}$. This does not depend on a possible
choice between $N|_{V_+}$ or $N|_{V_-}$ since both maps have the same
Jordan type.

Next, we consider \[\Dc = \bigcup_{r \in \NN^*} \{ (d_1, \dots, d_r)
\in \NN^r \mid d_1 \geq d_2 \geq \dots \geq d_r \geq 1 \} \] Define $d
: \Lambda \to \Dc$ by $d(\lambda) = (d_1(\lambda), \dots,
d_{r_\lambda}(\lambda))$. It is clear that $\Phi \circ d = m$, where
$\Phi : \Dc \to \NN$ is the map defined by $\Phi(d_1, \dots, d_r) =
\sum_{i=1}^r d_i$.

Finally, we can associate to $C \in \ok(n)$ a triple $(\Lambda, m, d)$
defined as above.

\begin{defn}\label{6.9bis}

  Let $\Jc_p$ be the set of all triples $(\Lambda, m, d)$ such that:

\begin{enumerate}

\item $\Lambda$ is a subset of $\CC \setminus \{0\}$ with
  $\sharp \Lambda \leq 2p$ and $\lambda \in \Lambda$ if, and only if,
  $-\lambda \in \Lambda$.

\item $m : \Lambda \to \NN^*$ satisfies $m(\lambda) = m(-\lambda)$,
  for all $\lambda \in \Lambda$ and $\sum_{\alpha \in \Lambda}
  m(\lambda) = 2p$.

\item $d : \Lambda \to \Dc$ satisfies $d(\lambda) = d(-\lambda)$, for
  all $\lambda \in \Lambda$ and $\Phi \circ d = m$.

\end{enumerate}

\end{defn}

Let $\Ic(2p)$ be the set of invertible elements in $\ok(2p)$ and
$\Ich(2p)$ be the set of $\OO(2p)$-adjoint orbits of elements in
$\Ic(2p)$. By the preceding remarks, there is a map $i : \Ic(2p) \to
\Jc_p$. The following Proposition classifies $\Ich(2p)$:

\begin{prop} \hfill \label{6.8} 

  The map $i : \Ic (2p) \to \Jc_p$ induces a bijection $\tilde{i} :
  \Ich(2p) \to \Jc_p$.
\end{prop}

\begin{proof}
  Let $C$ and $C' \in \Ic(2p)$ such that $C' = U \ C \ U^{-1}$ with $U
  \in \OO(2p)$. Let $S$, $S'$, $N$, $N'$ be respectively the
  semi-simple and nilpotent parts of $C$ and $C'$. Write $i(C) =
  (\Lambda, m, \lambda)$ and $i(C') = (\Lambda', m', \lambda')$.

  Then $S' = U \ S \ U^{-1}$. So $\Lambda' = \Lambda$ and $m' =
  m$. Also, $U(V_\lambda) = V'_\lambda$, for all $\lambda \in
  \Lambda$. Since $N' = U \ N \ U^{-1}$, then $N'|_{V'_\lambda} =
  U|_{V_\lambda} \ N|_{V_\lambda} \ U^{-1}|_{V_\lambda}$. Hence,
  $N|_{V_\lambda}$ and $N'|_{V_\lambda}$ have the same Jordan
  decomposition, so $d = d'$ and $\tilde{i}$ is well defined.

  To prove that $\tilde{i}$ is onto, we start with $\Lambda = \{
  \lambda_1, -\lambda_1, \dots, \lambda_k, -\lambda_k \}$, $m$ and $d$
  as in Definition \ref{6.9bis}. Define on the canonical basis: \[ S =
  \diag_{2p} (\overbrace{\lambda_1, \dots, \lambda_1}^{m(\lambda_1)},
  \dots, \overbrace{\lambda_k, \dots, \lambda_k}^{m(\lambda_k)},
  \overbrace{-\lambda_1, \dots, -\lambda_1}^{m(\lambda_1)}, \dots,
  \overbrace{-\lambda_k, \dots, -\lambda_k}^{m(\lambda_k)}). \] For
  all $1 \leq i \leq k$, let $d(\lambda_i) = (d_1(\lambda_i) \geq
  \dots d_{r_{\lambda_i}}(\lambda_i) \geq 1)$ and define \[
  N_+(\lambda) = \diag_{d(\lambda_i)} \left( J_{d_1(\lambda_i)},
    J_{d_2(\lambda_i)}, \dots, J_{d_{r_{\lambda_i}}(\lambda_i)}
  \right) \] on the eigenspace $V_{\lambda_i}$ and $0$ on the
  eigenspace $V_{-\lambda_i}$ where $J_d$ is the Jordan block of size
  $d$.

  By Lemma \ref{6.6}, $N(\lambda_i) := N_+(\lambda_i) -
  N_+(\lambda_i)^*$ is skew-symmetric on $V_{\lambda_i} \oplus
  V_{-\lambda_i}$. Finally, \[ \CC^{2p}= \oplusp_{i=1}^k \left(
    V_{\lambda_i} \oplus V_{-\lambda_i} \right).\] Define $N \in
  \ok(2p)$ by $N\left(\sum_{i=1}^k v_i \right) = \sum_{i=1}^k
  N(\lambda_i) (v_i)$, $v_i \in V_{\lambda_i} \oplus V_{-\lambda_i}$
  and $C = S + N \in \ok(2p)$. By construction, $i(C) = (\Lambda,
  m,d)$, so $\tilde{i}$ is onto.

  To prove that $\tilde{i}$ is one-to-one, assume that $C$, $C' \in
  \Ic(2p)$ and that $i(C) = i(C') = (\Lambda, m, d)$. Using the
  previous notation, since their respective semi-simple parts $S$ and
  $S'$ have the same spectrum and same multiplicities, there exist $U
  \in \OO(2p)$ such that $S' = U S U^{-1}$. For $\lambda \in \Lambda$,
  we have $U(V_\lambda) = V'_\lambda$ for eigenspaces $V_\lambda$ and
  $V'_\lambda$ of $S$ and $S'$.

  Also, for $\lambda \in \Lambda$, if $N$ and $N'$ are the nilpotent
  parts of $C$ and $C'$, then $N''(V_\lambda) \subset V_\lambda$, with
  $N'' = U^{-1} N' U$. Since $i(C) = i(C')$, then $N|_{V_\lambda}$ and
  $N'|_{V'_\lambda}$ have the same Jordan type. Since $ N'' = U^{-1}
  N' U$, then $N''|_{V_\lambda}$ and $N'|_{V'_\lambda}$ have the same
  Jordan type. So $N|_{V_\lambda}$ and $N''|_{V_\lambda}$ have the
  same Jordan type. Therefore, there exists $D_+ \in \Lc(V_\lambda)$
  such that $N''|_{V_\lambda} = D_+ N|_{V_\lambda} D_+^{-1}$. By Lemma
  \ref{6.6}, there exists $D(\lambda) \in \OO(V_\lambda \oplus
  V_{-\lambda})$ such that \[ N''|_{V_\lambda \oplus V_{-\lambda}} =
  D_+(\lambda) N|_{V_\lambda \oplus V_{-\lambda}} D_+(\lambda)^{-1}.\]
  We define $D \in \OO(2p)$ by $D|_{V_\lambda \oplus V_{-\lambda}} =
  D(\lambda)$, for all $\lambda \in \Lambda$. Then $N'' = D N D^{-1}$
  and $D$ commutes with $S$. Then $S' = (UD) S (UD)^{-1}$ and $N' =
  (UD) N (UD)^{-1}$ and we conclude \[ C' = (UD) C (UD)^{-1}.\]

\end{proof}

The classification of $\Sih(2p+2)$ can be deduced from the
classification of the set of orbits $\Ich(2p)$ by $\Jc_p$ as follows:
introduce an action of the multiplicative group $\CC^* = \CC \setminus
\{0\}$ on $\Jc_p$ by \[ \text{for all } \mu \in \CC^*, \ \mu \cdot
(\Lambda, m, d) = ( \mu \Lambda, m',d'), \ \forall \ (\Lambda, m, d)
\in \Jc_p, \lambda \in \Lambda, \] where $m'(\mu \lambda) = m
(\lambda), d'(\mu \lambda) = d(\lambda), \ \forall \ \lambda \in
\Lambda.$ We have $i(\mu C) = \mu i(C)$, for all $C \in \Ic(2p)$ and
$\mu \in \CC^*$. Hence, there is a bijection $\hat{i} : \ps(\Ich(2p))
\to \Jc_p / \CC^*$ given by $\hat{i} ([C]) = [i(C)]$, if $[C]$ is the
class of $C \in \Ic(2p)$ and $[(\Lambda, m, d)]$ is the class of
$(\Lambda, m, d) \in \Jc_p$.

\begin{prop}\label{6.11}
  The set $\Sih(2p+2)$ is in bijection with $\Jc_p / \CC^*$.
\end{prop}

\begin{proof}
  By Proposition \ref{4.11}, there is a bijection between
  $\Ssoi(2p+2)$ and $\widetilde{\ps(\ok(2p))}$. By restriction, that
  induces a bijection between $\Sih^{\mathrm{i}}(2p+2)$ and
  $\widetilde{\ps(\Ic(2p))}$. By Lemma \ref{5.7}, we have
  $\Sih^{\mathrm{i}}(2p+2) = \Sih(2p+2)$. Then, the result follows:
  given $\g \in \Si(2p+2)$ and an associated $\cb \in \Ic(2p)$, the
  bijection maps $\widetilde{\g}$ to $[i(\cb)]$ where $\widetilde{\g}$
  is the isomorphism class of $\g$.
\end{proof}

\begin{rem}
  Any $\g \in \Ss(n+2)$ can be decomposed as an amalgamated product of
  its Fitting components, $\g = \g_N \ta \g_I$ (Remark
  \ref{6.2}). Also, $\g \simeq \g'$ if, and only if, $\g_N \simeq
  \g_N'$ and $\g_I \simeq \g_I'$. Remark that $\g_N \in \Ns(k+2)$ for
  some $k \leq n$ and $\g_I \in \Si(2 \ell +2)$ for some $\ell$ with
  $2 \ell \leq n$ and $k + 2 \ell = n$. Up to isomorphism (or the
  equivalent notion of i-isomorphism, see Proposition \ref{6.4}), the
  classification of $\Ns(k+2)$ is known (Proposition \ref{5.4}) and
  the classification of $\Si(2 \ell +2)$ is known as well (Proposition
  \ref{6.11}). The decomposition of $\g_N$ and $\g_I$ as amalgamated
  products of Jordan-type Lie algebras is obtained in Propositions
  \ref{5.4} and \ref{6.5} and that allows us to write explicitly the
  commutation rules of $\g$. So, the complete description and
  classification (up to isomorphism or i-isomorphism) of $\Sso(n+2)$
  is achieved.

  Remark that aside the singular quadratic Lie algebras context, we
  can completely solve the problem of the classification of
  $\OO(n)$-adjoint orbits in $\ok(n)$ as follows: for $C \in \ok(n)$,
  consider its Fitting components $C_N$ and $C_I$. They belong
  respectively to $\Nc(k)$, $k \leq n$ and to $\Ic(2 \ell)$, $\ell
  \leq n$ with $k + 2 \ell = n$. Moreover, $C$ and $C'$ are conjugate
  if, and only if, $C_N$, $C_N'$ and $C_I$, $ C_I'$ are conjugate (it
  results from the proof of Proposition \ref{6.4}). But $C_N$ is
  nilpotent and the classification of nilpotent orbits is known (see
  Lemma \ref{5.3}). For the invertible $C_I$, the classification is
  given in Proposition \ref{6.8}. A Jordan-type decomposition of $C$
  can be then deduced (see (\ref{5.2}) and the proof of Proposition
  \ref{6.5}). This gives an explicit description and classification of
  $\OO(n)$-adjoint orbits in $\ok(n)$.
\end{rem}

\section{Quadratic dimension of reduced singular quadratic Lie
  algebras and invariance of $\mathrm{dup}(\g)$}

\subsection{} Let $(\g,B)$ be a quadratic Lie algebra. It is shown in
\cite{BB} that the space of invariant symmetric bilinear forms on $\g$
and the space generated by non-degenerated ones are the same. Let us
call it $\Bs(\g)$. The dimension of $\Bs(\g)$ is the {\em quadratic
  dimension} of $\g$, denote it by $d_q(\g)$. Obviously, $d_q(\g) =1$
if $\g$ is simple. If $\g$ is reductive, but neither simple, nor
one-dimensional, then \[ d_q(\g) = s(\g) + \dfrac{\dim(\Zs(\g))
  (1+\dim(\Zs(\g))}{2}, \] where $\Zs(\g)$ is the center of $\g$ and
$s(\g)$ is the number of simple ideals of a Levi factor of $\g$
\cite{BB}. A general formula for $d_q(\g)$ is not known. Here, we give
a formula for reduced singular quadratic Lie algebras. To any
symmetric bilinear form $B'$ on $\g$, there is an associated symmetric
map $D :\g \to \g$ satisfying \[B'(X,Y) = B(D(X),Y), \forall \ X,Y \in
\g.\] The following Lemma is straightforward.

\begin{lem}
  Let $(\g,B)$ be a quadratic Lie algebra, $B'$ be a bilinear form on
  $\g$ and $D \in \Lc(\g)$ its associated symmetric map. Then:

\begin{enumerate}

\item $B'$ is invariant if, and only if, $D$ satisfies 
\begin{equation} \label{etoile}
  D([X,Y]) =
  [D(X), Y] = [X, D(Y)], \ \forall \ X,Y \in \g.
\end{equation}

\smallskip

\item $B'$ is non-degenerate if, and only if, $D$ is invertible.

\smallskip

\end{enumerate}

\end{lem}

A symmetric map $D$ satisfying \ref{etoile} is called a {\em
  centromorphism} of $\g$. The space of centromorphisms and the space
generated by invertible centromorphisms are the same, denote it by
$\Cs(\g)$. We have $d_q(\g) = \dim(\Cs(\g))$.

\begin{prop} \label{7.2} Let $\g$ be a reduced singular quadratic Lie
  algebra and $D \in \Lc(\g)$ be a symmetric map. Then:

\begin{enumerate}

\item $D$ is a centromorphism if, and only if, there exists $\mu \in
  \CC$ and a symmetric map $\Zb : \g \to \Zs(\g)$ such that
  $\Zb|_{[\g,\g]} = 0$ and $D = \mu \Id + \Zb$. Moreover $D$ is
  invertible if, and only if, $\mu \neq 0$.

\smallskip

\item \[ d_q(\g) = 1 + \dfrac{\dim(\Zs(\g)) (1+\dim(\Zs(\g))}{2}. \]

\smallskip

\end{enumerate}

\end{prop}

\begin{proof} \hfill

\begin{enumerate}

\item If $\g = \ok(3)$, with $B = \lambda \kappa$ and $\kappa$ the
  Killing form, the two results are obvious. So, we examine the case
  where $\g$ is solvable, and then $\g$ can be realized as a double
  extension: $\g = (\CC X_1 \oplus \CC Y_1) \oplusp \qk$, with
  corresponding bilinear form $\overline{B}$ on $\qk$, $C = \ad(Y_1)$,
  $\cb = C|_\qk \in \ok(\qk)$.

  Let $D$ be an invertible centromorphism. One has $D \circ \ad(X) =
  \ad(X) \circ D$, for all $X \in \g$ and that implies $D C = C
  D$. Using formula (1) of Lemma \ref{4.3} and $CD = DC$, from $[D(X),
  Y_1] = [X, D(Y_1)]$, we find $D(C(X)) = $ \linebreak $B(D(X_1), Y_1)
  C(X)$. Let $\mu = B(D(X_1), Y_1)$. Since $D$ is invertible, one has
  $\mu \neq 0$ and $C (D - \mu \Id) = 0$. Since $\ker(C) = \CC X_1
  \oplus \ker(\cb) \oplus \CC Y_1 = \Zs(\g) \oplus \CC Y_1$, there
  exists a map $\Zb : \g \to \Zs(\g)$ and $\varphi \in \g^*$ such that
  $D - \mu \Id = \Zb + \varphi \otimes Y_1$. But $D$ maps $[\g,\g]$
  into itself, so $\varphi|_{[\g,\g]} = 0$. One has $[\g,\g] = \CC X_1
  \oplus \im(\cb)$. If $X \in \im(\cb)$, let $X = C(Y)$. Then $D(X) =
  D(C(Y)) = \mu C(Y)$, so $D(X) = \mu X$. For $Y_1$, $D([Y_1,X]) =
  DC(X) = \mu C(X)$ for all $X \in \g$. But also, $D([Y_1,X]) =
  [D(Y_1),X] = \mu C(X) + \varphi(Y_1) C(X)$, hence $\varphi(Y_1) =
  0$.

  Assume we have shown that $D(X_1) = \mu X_1$. Then if $X \in \qk$,
  $B(D(X_1), X) = \mu B(X_1,X) = 0$. Moreover, $B(D(X_1), X) = B(X_1,
  D(X))$, so $\varphi(X) = 0$. Thus, to prove (1), we must prove that
  $D(X_1) = \mu X_1$. We decompose $\qk$ respectively to $\cb$ as in
  Appendix 1. Let $\lk = \ker(\cb)$. Then: \[\qk = (\lk \oplus \lk')
  \oplusp (\uk \oplus \uk')\] and $C$ is an isomorphism from $\lk'
  \oplusp (\uk \oplus \uk')$ onto $\lk \oplusp (\uk \oplus
  \uk')$. Or \[\qk = (\lk + \lk') \oplusp \CC T \oplusp (\uk \oplus
  \uk')\] and $C$ is an isomorphism from $\lk' \oplusp \CC T \oplusp
  (\uk \oplus \uk')$ onto $\lk \oplusp \CC T \oplusp (\uk \oplus
  \uk')$.

  If $\uk \oplus \uk' \neq \{0\}$, there exist $X'$, $Y' \in \uk
  \oplus \uk'$ such that $B(X',Y') = -1$ and $X$, $Y \in \lk' \oplusp
  (\uk \oplus \uk')$ (resp. $\lk' \oplusp \CC T \oplusp (\uk \oplus
  \uk')$) such that $X' = C(X)$, $Y' = C(Y)$. It follows that $[C(X),
  Y] = X_1$ and then $D(X_1) = [D C(X), Y] = \mu [C(X),Y] = \mu X_1$.

  If $\uk \oplus \uk' = \{0\}$, then either $\qk = (\lk + \lk')
  \oplusp \CC T$ or $\qk = \lk + \lk'$. The first case is similar to
  the situation above, setting $X' = Y' = \dfrac{T}{i}$ and $X$, $Y
  \in \lk' \oplusp \CC T$. In the second case, $\lk = \im(\cb)$ is
  totally isotropic and $C$ is an isomorphism from $\lk'$ onto
  $\lk$. For any non-zero $X \in \lk'$, choose a non-zero $Y \in \lk'$
  such that $B(C(X),Y) = 0$. Then $D([X,Y]) = D(B(C(X),Y) X_1) =
  0$. But this is also equal to $[D(X), Y] = \mu [X,Y] + \varphi(X)
  C(Y)$. Since $D$ is invertible, $[X,Y] = 0$ and we conclude that
  $\varphi(X) = 0$. Therefore $\varphi|_{\lk'} = 0$. There exist $L$,
  $L' \in \lk'$ such that $X_1 = [L, L']$ and then $D(X_1) = \mu X_1$.

  Finally, $\Cs(\g)$ is generated by invertible centromorphism, so the
  necessary condition of (1) follows. The sufficiency is a simple
  verification.

 \smallskip

\item As in (1), we can restrict to a double extension and follow the
  same notation. By (1), $D$ is a centromorphism if, and only if,
  $D(X) = \mu X + \Zb(X)$, for all $X \in \g$ with $\mu \in \CC$ and
  $\Zb$ is a symmetric map from $\g$ into $\Zs(\g)$ satisfying
  $\Zb|_{[\g,\g]} = 0$. To compute $d_q(\g)$, we use Appendix
  1. Assume $\dim(\qk)$ is even and write $\qk = (\lk \oplus \lk')
  \oplusp (\uk \oplus \uk')$ with $\lk = \ker(\cb)$, $\Zs(\g) = \CC
  X_1 \oplus \lk$, $\im(\cb) = \lk \oplusp (\uk \oplus \uk')$ and
  $[\g,\g] = \CC X_1 \oplus \im(\cb)$. Let us define $\Zb : \lk'
  \oplusp \CC Y_1 \to \lk \oplusp \CC X_1$: set basis $\{X_1, X_2,
  \dots, X_r \}$ of $\lk \oplus \CC X_1$ and $\{Y_1' = Y_1, Y_2',
  \dots, Y_r' \}$ of $\lk' \oplus \CC Y_1$ such that $B(Y_i', X_j) =
  \delta_{ij}$. Then $\Zb$ is completely defined by \[\Zb \left(
    \sum_{j=1}^r \mu_j Y_j' \right) = \sum_{i=1}^r \left( \sum_{j=1}^r
    \nu_{ij} \mu_j \right) X_i\] with $\nu_{ij} = \nu_{ji} = B(Y_i',
  \Zb(Y_j'))$ and the formula follows. The case of $\dim(\qk)$ odd is
  completely similar.

\smallskip

\end{enumerate}

\end{proof}

\subsection{} As a consequence of Proposition \ref{7.2}, we prove:

\begin{prop} \label{7.3}
  The $\dup$-number is invariant under isomorphism, i.e. if $\g$ and
  $\g'$ are quadratic Lie algebras with $\g \simeq \g'$, then
  $\dup(\g) = \dup(\g')$.
\end{prop}

\begin{proof}
  Assume that $\g \simeq \g'$. Since an i-isomorphism does not change
  $\dup(\g')$, we can assume that $\g = \g'$ as Lie algebras equipped
  with invariant bilinear forms $B$ and $B'$. Thus, we have two
  $\dup$-numbers, $\dup_B(\g)$ and $\dup_{B'}(\g)$.

  We choose $\zk$ such that $\Zs(\g) = \left( \Zs(\g) \cap [\g,\g]
  \right) \oplus \zk$. Then $\zk \cap \zk^{\perp_B} = \{0\}$, $\zk$ is
  a central ideal of $\g$ and $\g = \lk \opluspb \zk$ with $\lk$ a
  reduced quadratic Lie algebra. Then $\dup_B(\g) = \dup_B(\lk) $ (see
  (\ref{2.2s})). Similarly, $\zk \cap \zk^{\perp_{B'}} = \{0\}$, $\g =
  \lk' \opluspbp \zk$ with $\lk$ a reduced quadratic Lie algebra and
  $\dup_{B'}(\g) = \dup_{B'}(\lk')$. Now, $\lk$ and $\lk'$ are
  isomorphic to $\g / \zk$, so $\lk \simeq \lk'$. Therefore, it is
  enough to prove the result for reduced quadratic Lie algebras to
  conclude that $\dup_B(\lk) = \dup_{B'}(\lk)$ and then that
  $\dup_B(\g) = \dup_{B'}(\g)$.

  Consider $\g$ a reduced quadratic Lie algebra equipped with bilinear
  forms $B$ and $B'$ and associated 3-forms $I$ and $I'$. (see
  (\ref{1.5})). We have $\dup_B(\g) = \dim(\Vs_I)$ and $\dup_{B'}(\g)
  = \dim(\Vs_{I'})$ with $\Vs_I = \{ \alpha \in \g^* \mid \alpha
  \wedge I = 0 \}$ and $\Vs_{I'} = \{ \alpha \in \g^* \mid \alpha
  \wedge I' = 0 \}$.

  We start with the case $\dup_B(\g) =3$. This is true if, and only
  if, $\dim([\g,\g]) = 3$ \cite{PU07}. Then $\dup_{B'} (\g) = 3$.

  If $\dup_B(\g) = 1$, then $\g$ is of type $\Sb_1$ with respect to
  $B$. We apply Proposition \ref{7.2} to obtain an invertible
  centromorphism $D = \mu \Id + \Zb$ for a non-zero $\mu \in \CC$, $\Zb =
  \g \to \Zs(\g)$ satisfying $\Zb|_{[\g,\g]} = 0$ and such that $B'(X,Y)
  = B(D(X), Y)$, for all $X, Y \in \g$. Then $I'(X,Y,Z) = B'([X,Y],Z)
  = B([D(X),Y],Z) = \mu B([X,Y],Z) = \mu I(X,Y,Z)$, for all $X$, $Y$,
  $Z \in \g$. So $I' = \mu I$ and $\dup_{B'}(\g) = \dup_B(\g)$. 

  Finally, if $\dup_B(\g) = 0$, then from the previous cases, $\g$
  cannot be of type $\Sb_3$ or $\Sb_1$ with respect to $B'$, so
  $\dup_{B'}(\g) = 0$.

\end{proof}

\section{Appendix 1}

In this Appendix, we recall some facts on skew-symmetric maps used in
the paper. Nothing here is new, but short proofs are given for the sake of
completeness.

Throughout this section, $(V,B)$ is a quadratic vector space and $C$
is an element of $\ok(V)$. We recall the useful identity $\ker(C) =
(\im(C))^\perp$.

\begin{lem} \label{A1} There exist subspaces $W$ and $N$ of $V$ such
  that:

\begin{enumerate}

\item $N \subset \ker(C)$, $C(W) \subset W$ and $V = W \oplusp N$.

\smallskip

\item Let $B_W = B|_{W \times W}$ and $C_W = C|_{W}$.  Then $B_W$ is
  non-degenerate, $C_W \in \ok(W, B_W)$ and $\ker(C_W) \subset \im(C_W)
  = \im(C)$.

\end{enumerate}

\end{lem}

\begin{proof}
  We follow the proof of Proposition \ref{2.8}, given in
  \cite{PU07}. Let $N_0 = \ker(C) \cap \im(C)$ and let $N$ be a
  complementary subspace of $N_0$ in $\ker(C)$, $\ker(C) = N_0 \oplus
  N$. Since $\ker(C) = (\im(C))^\perp$, we have $B(N_0,N) =\{ 0 \}$
  and $N \cap N^\perp = \{ 0 \}$. So, if $W = N^\perp$, one has $V = W
  \oplusp N$. From $C(N) = \{0\}$, we deduce that $C(W) \subset W$.

  It is clear that $B$ is non-degenerate and that $C_W \in
  \ok(W)$. Moreover, since $C(W) \subset W$ and $C(N) = \{0\}$, then
  $\im(C) = \im(C_W)$. It is immediate that $\ker(C_W) = N_0$, so
  $\ker(C_W) \subset \im(C_W)$.
\end{proof}

\begin{lem} Assume that $\ker(C) \subset \im(C)$. Denote $L =
  \ker(C)$. Let $\{L_1, \dots, L_r \}$ be a basis of $L$.

\begin{enumerate}

\item If $\dim(V)$ is even, there exist subspaces $L'$ with basis
  $\{L_1', \dots, L_r' \}$, $U$ with basis $\{U_1, \dots, U_s \}$ and
  $U'$ with basis $\{U_1', \dots, U_s' \}$ such that $B(L_i, L_j') =
  \delta_{ij}$, for all $1 \leq i,j \leq r$, $L$ and $L'$ are totally
  isotropic, $B(U_i, U_j') = \delta_{ij}$, for all $1 \leq i,j \leq
  s$, $U$ and $U'$ are totally isotropic and \[ V= (L \oplus L')
  \oplusp (U \oplus U').\] Moreover $\im(C) = L \oplusp (U \oplus U')$
  and $C : L' \oplusp (U \oplus U') \to L \oplusp (U \oplus U')$ is a
  bijection.

\smallskip

\item If $\dim(V)$ is odd, there exist subspaces $L'$, $U$ and $U'$ as
  in (1) and $v \in V$ such that $B(v,v) = 1$ and \[ V= (L \oplus L')
  \oplusp \CC v \oplusp (U \oplus U').\] Moreover $\im(C) = L \oplusp
  \CC v \oplusp (U \oplus U')$ and $C : L' \oplusp \CC v \oplusp (U
  \oplus U') \to L \oplusp \CC v \oplusp (U \oplus U')$ is a
  bijection.

\smallskip

\item In both cases, $\rank(C)$ is even.

\end{enumerate}

\end{lem}

\begin{proof}

  Since $\left( \ker(C)\right)^\perp = \im(C)$, $L$ is isotropic.

\begin{enumerate}

\item If $\dim(V)$ is even, there exist maximal isotropic subspaces
  $W_1$ and $W_2$ such that $V = W_1 \oplus W_2$ \cite{Bour59} and $L
  \subset W_1$. Let $U$ be a complementary subspace of $L$ in $W_1$,
  $W_1 = L \oplus U$ and $\{U_1, \dots, U_s \}$ a basis of
  $U$. Consider the isomorphism $\Psi : W_2 \to W_1^*$ defined by
  $\Psi(w_2) (w_1) = B(w_2, w_1)$, for all $w_1 \in W_1$, $w_2 \in
  W_2$. Define $L_i' = \psi^{-1}(L_i^*)$, $1 \leq i \leq r$, $L' =
  \spa \{ L_1', \dots, L_r' \}$, $U_j' = \psi^{-1}(U_j^*)$, $1 \leq j
  \leq s$, $U' = \spa \{ U_1', \dots, U_s' \}$. Then $B(L_i, L_j') =
  \delta_{ij}$, $1 \leq i,j \leq r$, $L$ and $L'$ are isotropic,
  $B(U_i, U_j') = \delta_{ij}$, for all $1 \leq i,j \leq s$, $U$ and
  $U'$ are isotropic and \[ V= (L \oplus L') \oplusp (U \oplus U').\]
  Since $\im(C) = L^\perp$, we have $\im(C) = L \oplusp (U \oplus
  U')$. Finally, if \linebreak $v \in L' \oplusp (U \oplus U')$ and
  $C(v) =0$, then $v \in L$. So $v =0$. Therefore $C$ is one to one
  from $L' \oplusp (U \oplus U')$ into $L \oplusp (U \oplus U')$ and
  since the dimensions are the same, $C$ is a bijection.

\smallskip

\item There exist maximal isotropic subspaces $W_1$ and $W_2$ such
  that $V = (W_1 \oplus W_2) \oplusp \CC v$, with $v \in V$ such that
  $B(v,v)=1$ and $L \subset W_1$ \cite{Bour59}. Then the proof is
  essentially the same as in (1).

\item Assume $\dim(V)$ even. Define a bilinear form $\Delta$ on $L'
  \oplusp (U \oplus U')$ by $\Delta(v_1, v_2) = B(v_1, C(v_2))$, for
  all $v_1$, $v_2 \in L' \oplusp (U \oplus U')$. Since $C \in \ok(V)$,
  $\Delta$ is skew-symmetric. Let $v_1 \in L' \oplusp (U \oplus U')$
  such that $\Delta( v_1, v_2 ) = 0$, for all $v_2 \in L' \oplusp (U
  \oplus U')$. Then $B(v_1, w)= 0$, for all $w \in L \oplusp (U \oplus
  U')$. It follows that $B(v_1, w) = 0$, for all $w \in V$, so $v_1 =
  0$ and $\Delta$ is non-degenerate. So $\dim(L' \oplusp (U \oplus
  U')$ is even. Therefore $\dim(L') = \dim(L)$ is even and $\rank(C)$
  is even. If $V$ is odd-dimensional, the proof is completely similar.

\end{enumerate}
\end{proof}

\begin{cor}
  If $C \in \ok(V)$, then $\rank(C)$ is even.
\end{cor}

\begin{proof}
  By Lemma \ref{A1}, $\im(C) = \im(C_W)$ and $\rank(C_W)$ is even by
  the preceding Lemma. 
\end{proof}

For instance, if $C \in \ok(V)$ and $C$ is invertible, then $\dim(V)$
must be even. But this can also be proved directly: when $C$ is
invertible, then the skew-symmetric form $\Delta_C$ on $V$ defined by
$\Delta_C(v_1, v_2) = B(v_1, C(v_2))$, for all $v_1$, $v_2 \in V$, is
clearly non-degenerate.

When $C$ is semi-simple (i.e. diagonalizable), we have $V = \ker(C)
\oplusp \im(C)$ and $C|_{\im(C)}$ is invertible. So semi-simple
elements are completely described by:

\begin{lem} \label{A3} Assume $C$ is semi-simple and invertible. Then
  there is a basis \linebreak $\{e_1, \dots,e_p, f_1, \dots,f_p\}$ of
  $V$ such that $B(e_i,e_j) = B(f_i,f_j) = 0$, $B(e_i,f_j) =
  \delta_{ij}$, $1 \leq i,j \leq p$. For $1 \leq i \leq p$, there
  exist non-zero $\lambda_i \in \CC$ such that $C(e_i) = \lambda_i
  e_i$ and $C(f_i) = - \lambda_i f_i$.

  Moreover, if $\Lambda$ denotes the spectrum of $C$, then $\lambda
  \in \Lambda$ if, and only if, $-\lambda \in \Lambda$; $\lambda$ and
  $-\lambda$ have the same multiplicity.
\end{lem}

\begin{proof}
  We prove the result by induction on $\dim(V)$. Assume $\dim(V) =
  2$. Let $\{e_1, e_2 \}$ be an eigenvector basis of $V$ corresponding
  to eigenvalues $\lambda_1$ and $\lambda_2$. We have $B(C(v), v') = -
  B(v,C(v'))$ and $C$ is invertible, so $B(e_1, e_1) = B(e_2, e_2) =
  0$, $B(e_1,e_2) \neq 0$ and $\lambda_2 = - \lambda_1$. Let $f_1 =
  \dfrac{1}{B(e_1,e_2)} e_2$, then the basis $\{e_1, f_1 \}$ is a
  convenient basis.

  Assume that the result is true for quadratic vector spaces of
  dimension $n$ with $n \leq 2(p-1)$. Assume $\dim(V) = 2p$. Let
  $\{e_1, \dots, e_{2p}\}$ be an eigenvector basis with corresponding
  eigenvalues $\lambda_1, \dots, \lambda_{2p}$. As before, $B(e_i, e_i
  ) = 0$, $1 \leq i \leq 2p$, so there exists $j$ such that $B(e_1,
  e_j) \neq 0$. Then $\lambda_j = - \lambda_1$. Let $f_1 =
  \dfrac{1}{B(e_1,e_j)} e_j$. Then $B|_{\spa \{e_1,f_1\}}$ is
  non-degenerate, so $V = \spa \{e_1, f_1 \} \oplusp V_1$, where $V_1
  = \spa \{ e_1, f_1 \}^\perp$. But $C$ maps $V_1$ into itself, so we
  can apply the induction assumption and the result follows.
\end{proof}

As a consequence, we have this classical result, used in Section 5:

\begin{lem}\label{A4} \hfill

\begin{enumerate}

\item Let $C$ be a semi-simple element of $\ok(n)$. Then $C$ belongs
  to the $\SO(n)$-adjoint orbit of an element of the standard Cartan
  subalgebra of $\ok(n)$ (i.e., an element with matrix $\diag_{2p}
  (\lambda_1, \dots, \lambda_p, -\lambda_1, \dots, -\lambda_p)$ if $n
  =2p$ and $\diag_{2p+1} (\lambda_1, \dots, \lambda_p,0,- \lambda_1,
  \dots, -\lambda_p)$ if $n =2p+1$ in the canonical basis of $\CC^n$).

\smallskip

\item Let $C$ and $C'$ be semi-simple elements of $\ok(n)$. Then $C$
  and $C'$ are in the same $\OO(n)$-adjoint orbit if, and only if,
  they have the same spectrum, with same multiplicities.

\end{enumerate}

\end{lem}

\begin{proof} \hfill

\begin{enumerate}

\item We have $\CC^n = \ker(C) \oplusp \im(C)$ and $\rank(C)$ is
  even. So $\dim(\ker(C))$ is even if $n = 2p$ and odd, if $n =
  2p+1$. Then apply Lemma \ref{A3} to $C|_{\im(C)}$ to obtain the
  result.

\smallskip

\item If $C$ and $C'$ have the same spectrum and their eigenvalues,
  same multiplicities, they are $\OO(n)$-conjugate to the same element
  of the standard Cartan subalgebra.

\end{enumerate}

\end{proof}

\begin{rem} \hfill

\begin{enumerate}

\item Attention: $\OO(n)$-adjoint orbits are generally not the same as
  $\SO(n)$-adjoint orbits. 

\smallskip

\item Lemma \ref{A4}(1) is a particular case of a general and
  classical result on semi-simple Lie algebras: any semi-simple
  element of a semi-simple Lie algebra belongs to a Cartan subalgebra
  and all Cartan subalgebras are conjugate under the adjoint action
  \cite{Sam}. Here, $\ok(n)$ is a semi-simple Lie algebra and the
  adjoint group is $\SO(n)$.

\end{enumerate}

\end{rem}

\section{Appendix 2}

Here we prove:

\begin{lem}
  Let $(\g,B)$ be a non-Abelian 5-dimensional quadratic Lie
  algebra. Then $\g$ is a singular quadratic Lie algebra.
\end{lem}

\begin{proof} \hfill

\begin{itemize}

\item We assume $\g$ is not solvable and we write $\g = \sk \oplus
  \rk$ with $\sk$ semi-simple and $\rk$ the radical of $\g$
  \cite{Bourg}. Then $\sk \simeq \slk(2)$ and $B|_{\sk \times \sk} =
  \lambda \kappa$ where $\kappa$ is the Killing form.

  If $\lambda = 0$, consider $\Psi : \sk \to \rk^*$ defined by
  $\Psi(S)(R) = B(S,R)$, for all $S \in \sk$, $R \in \rk$. Then $\Psi$
  is one-to-one and $\Psi\left( \ad(X)(S) \right) = \check{\ad}(X)
  (\psi(S))$, for all $X$, $S \in \sk$. So $\Psi$ must be a
  homomorphism from the representation $(\sk, \ad|_\sk)$ of $\sk$ into
  the representation $(\rk^*, \check{\ad}|_\sk)$, so $\Psi = 0$, a
  contradiction.

  So $\lambda \neq 0$. Then $B|_{\sk \times \sk}$ is
  non-degenerate. Therefore $\g = \sk \oplusp \sk^\perp$ and
  $\ad(\sk)|_{\sk^\perp}$ is an orthogonal 2-dimensional
  representation of $\sk$. Hence, \linebreak $\ad(\sk)|_{\sk^\perp}=0$
  and $[\sk, \sk^\perp] = 0$. We have $B(X,[Y,Z]) = B([X,Y], Z) = 0$,
  for all $X \in \sk$, $Y \in \sk^\perp$, $Z \in \g$. It follows that
  $\sk^\perp$ is an ideal of $\g$ and therefore a quadratic
  2-dimensional Lie algebra. So $\sk^\perp$ is Abelian. Finally, $\g =
  \sk \oplusp \sk^\perp$ with $\sk^\perp$ a central ideal of $\g$, so
  $\dup(\g) = \dup(\sk) = 3$.

\smallskip

\item We assume that $\g$ is solvable and we write $\g = \lk \oplusp
  \zk$ with $\zk$ a central ideal of $\g$ (Proposition
  \ref{2.8}). Then $\dim(\lk) \geq 3$. If $\dim(\lk) = 3$ or 4, then
  it is proved in Proposition \ref{2.6a} that $\lk$ is singular, so
  $\g$ is singular. So we can assume that $\g$ is reduced,
  i.e. $\Zs(\g) \subset [\g,\g]$. It results that $\dim(\Zs(\g)) = 1$
  or 2 (Remark \ref{2.3}).

\begin{itemize}

\item If $\dim(\Zs(\g)) = 1$, $\Zs(\g) = \CC X_0$. Then $\dim([\g,
  \g]) = 4$ and $[\g,\g] = X_0^\perp$. We can choose $Y_0$ such that
  $B(X_0, Y_0) =1$ and $B(Y_0, Y_0) = 0$. Let $\qk = (\CC X_0 \oplus
  \CC Y_0)^\perp$. Then $\g = (\CC X_0 \oplus \CC Y_0) \oplusp
  \qk$. If $X$, $X' \in \qk$, then $B(X_0, [X,X']) = B([X_0, X], X') =
  0$, so $[X,X'] \in X_0^\perp$. Write $[X,X'] = \lambda(X,X') X_0 +
  [X,X']_\qk$ with $[X, X']_\qk \in \qk$. Remark that $[X, [X', X'']]
  = \lambda(X,[X',X'']_\qk) X_0 + [X, [X',X'']_\qk]_\qk$, for all $X$,
  $X'$, $X'' \in \qk$. So $[\cdot, \cdot]_\qk$ satisfies the Jacobi
  identity. Moreover $B([X,X'], X'') = - B(X', [X, X'']_\qk)$. But
  also $B([X,X'], X'') = B([X,X']_\qk, X'')$. So $(\qk,[\cdot,
  \cdot]_\qk, B|_{\qk \times \qk})$ is a 3-dimensional quadratic Lie
  algebra.

  If $\qk$ is an Abelian Lie algebra, then $[X, X'] \in \CC X_0$, for
  all $X$, $X' \in \qk$. Write $B(Y_0, [X,X']) = B([Y_0, X], X')$ to
  obtain $[X,X'] = \linebreak B(\ad(Y_0)(X), X') X_0$, for all $X$,
  $X' \in \qk$. Since $\dim(\qk) = 3$ and $\ad(Y_0)|_\qk$ is
  skew-symmetric, there exists $Q_0 \in \qk$ such that $\ad(Y_0)(Q_0)
  = 0$. It follows that $Q_0 \in \Zs(\g)$ and that is a contradiction
  since $\dim(\Zs(\g)) = 1$.

  Therefore $(\qk, [\cdot, \cdot]_\qk) \simeq \slk(2)$.  Consider \[ 0
  \to \CC X_0 \to X_0^\perp \to \qk \to 0.\] Then there is a section
  $\sigma : \qk \to X_0^\perp$ such that \linebreak
  $\sigma([X,X']_\qk) = [\sigma(X), \sigma(X')]$, for all $X$, $X' \in
  \qk$ \cite{Bourg}. Then $\sigma(\qk)$ is a Lie subalgebra of $\g$,
  isomorphic to $\slk(2)$ and that is a contradiction since $\g$ is
  solvable.

\smallskip

\item If $\dim(\Zs(\g)) = 2$, then we choose a non-zero $X_0 \in
  \Zs(\g)$ and $Y_0 \in \g$ such that $B(X_0, Y_0) = 1$ and $B(Y_0,
  Y_0) = 0$. Let $\qk = (\CC X_0 \oplus \CC Y_0)^\perp$. Then $\g =
  (\CC X_0 \oplus \CC Y_0) \oplusp \qk$ and as in the preceding case,
  $[X,X'] \in X_0^\perp$, for all $X$, $X' \in \qk$.  Write $[X,X'] =
  \lambda(X,X') X_0 + [X,X']_\qk$ with $[X, X']_\qk \in \qk$.  Same
  arguments as in the preceding case allow us to conclude that
  $[\cdot, \cdot]_\qk$ satisfies the Jacobi identity and that $B|_{\qk
    \times \qk}$ is invariant. So $(\qk,[\cdot, \cdot]_\qk, B|_{\qk
    \times \qk})$ is a 3-dimensional quadratic Lie algebra.

  If $\qk \simeq \slk(2)$, then apply the same reasoning as in the
  preceding case to obtain a contradiction with $\g$ solvable.

  If $\qk$ is an Abelian Lie algebra, then $[X, X'] \in \CC X_0$, for
  all $X$, $X' \in \qk$. Again, as in the preceding case, $[X,X'] =
  B(\ad(Y_0)(X), X') X_0$, for all $X$, $X' \in \qk$. Then it is easy
  to check that $\g$ is a double extension of the quadratic vector
  space $\qk$ by $\cb = \ad(Y_0)|_\qk$. By Proposition \ref{4.4}, $\g$
  is singular.
\end{itemize}

\end{itemize}
\end{proof}

\begin{rem}

  Let us give a list of all non-Abelian 5-dimensional quadratic Lie
  algebras:

\begin{itemize}

\item $\g \iiso \ok(3) \oplusp \CC^2$ with $\CC^2$ central, $\ok(3)$
  equipped with bilinear form $\lambda \kappa$, $\lambda \in \CC$,
  $\lambda \neq 0$ and $\kappa$ the Killing form. We have $\dup(\g) =
  3$.

\smallskip

\item $\g \iiso \g_4 \oplusp \CC$ with $\CC$ central, $\g_4$ the
  double extension of $\CC$ by $\begin{pmatrix} 1 & 0 \\ 0 &
    -1 \end{pmatrix}$, $\g$ is solvable, non-nilpotent and $\dup(\g) =
  3$.

\smallskip

\item $\g \iiso \g_5$, a double extension of $\CC^3$ by
  $\begin{pmatrix} 0 & 1 & 0 \\ 0 & 0 & -1 \\ 0 & 0 &
    0 \end{pmatrix}$, $\g$ is nilpotent and $\dup(\g) = 3$.

\end{itemize}

\smallskip

See Proposition \ref{4.5} for the definition of $\g_4$ and
$\g_5$. Remark that $\g_4 \oplusp \CC$ is actually the double
extension of $\CC^3$ by $\begin{pmatrix} 1 & 0 & 0 \\ 0 & 0 & 0 \\ 0 &
  0 & -1 \end{pmatrix}$

\end{rem}

\bibliographystyle{amsxport}

\end{document}